%% file: ms.tex
\documentclass{amsart}
\usepackage{algtop}
\usepackage[all, cmtip]{xy}
\usepackage{caption}
\usepackage{graphics}
\usepackage{mathtools}
\usepackage{enumitem}
\usepackage[T1]{fontenc}
\usepackage[utf8]{inputenc}
\usepackage{makecell}
\usepackage{multirow}
\usepackage[dvipsnames,table]{xcolor}
\usepackage{url}
\usepackage[linewidth=1pt]{mdframed}

\usepackage{tikz}
\usetikzlibrary{arrows.meta}
\usetikzlibrary{calc}

\makeatletter
	\let\c@equation\c@theorem
\makeatother
\numberwithin{equation}{section}

\newcommand{\makebibliography}{
    \newpage
    \bibliographystyle{plain}
    \bibliography{ms.bib}
}

\newcommand{\AF}{\mathrm{AF}}

\newcommand{\ZESS}[1]{\prescript{#1\!}{}{Z}}
\newcommand{\BESS}[1]{\prescript{#1\!}{}{B}}

\newtheorem{proofstr}[theorem]{Proof Structure}

\definecolor{lightgray}{RGB}{200,200,200}
\definecolor{llgray}{RGB}{235,235,235}

\title[Machine Proofs for Adams Differentials and Extension Problems]{Machine Proofs for Adams Differentials and Extension Problems Among CW Spectra}
\author{Weinan Lin, Guozhen Wang and Zhouli Xu}
\date{}

\begin{document}
\begin{abstract}
In this document, we describe the process of obtaining numerous Adams differentials and extensions using computational methods, as well as how to interpret the dataset uploaded to Zenodo \cite{LWXZenodo}.
%\cite{}.%
Detailed proofs of the machine-generated results are also provided. The dataset includes information on 49 CW spectra, 180 maps, and 61 cofiber sequences. Leveraging these results, and with the addition of some ad hoc arguments derived through human insight, we successfully resolved the Last Kervaire Invariant Problem in dimension 126 \cite{LWX126}.
\end{abstract}

\maketitle
% \tableofcontents
% \blfootnote{The author is supported by China Postdoctoral Science Foundation 2021TQ0015 and the Fundamental Research Funds for the Central Universities, Peking University}

\section{Introduction}
In this project, we first compute the Adams $E_2$-pages and maps for many CW spectra using the computer program \texttt{./Adams}, developed by the first author \cite{LinProgram}. All relevant spectra and maps have been uploaded to the Zenodo repository \cite{LWXZenodo} and are described in detail in Section \ref{sec:E2}. Additionally, another program, \texttt{./ss}, also developed by the first author \cite{LinProgram}, computes Adams differentials and extensions. The resulting data are included in the Zenodo repository and described in Section \ref{sec:diff}.

The program \texttt{./ss} also generates proofs for all computed results, which are organized into a comprehensive resource we refer to as the Table of Proofs. This program partially employs the Generalized Leibniz Rule and the Generalized Mahowald Trick—two theorems established in the paper \cite{LWX126} addressing the Last Kervaire Invariant Problem. These theorems enable the program to derive new differentials and extensions from previously computed ones. Detailed instructions on how to interpret and navigate the Table of Proofs are provided in Section \ref{sec:proofs}.

In Section \ref{sec:charts}, we present several charts and tables of machine-generated Adams differentials that are particularly relevant to the solution of the Last Kervaire Invariant Problem.

\section{Adams $E_2$ data}\label{sec:E2}
\subsection{Adams $E_2$-pages}
Tables \ref{tb:CW1} provides a list of 49 CW spectra considered in this work, along with the maximum internal degrees computed for the Adams $E_2$-pages of the corresponding spectra.

\input{cw_data}\vspace{0.5cm}

For a CW spectrum \texttt{X}, the Adams $E_2$-page data for \texttt{X} is available in the database file \texttt{X\_AdamsSS.db} within \texttt{kervaire\_database.zip}. Alternatively, the data can be accessed as three plain-text table files:
\begin{itemize}
    \item \texttt{X\_AdamsE2\_generators.csv} (for generators),
    \item \texttt{X\_AdamsE2\_relations.csv} (for relations),
    \item \texttt{X\_AdamsE2\_basis.csv} (for $\bF_2$-basis).
\end{itemize}
These files are included in \texttt{kervaire\_csv.rar}, available in the Zenodo repository.
We illustrate these tables with examples. Tables \ref{tb:gen}--\ref{tb:basis} present the first few rows of the corresponding files for the spectrum \texttt{C2} (the naming convention will be explained later).

\begin{table}
\begin{tabular}{|c|c|c|c|}\hline
id & name & stem & s \\\hline
0 & \verb|[0]| & 0 & 0 \\\hline
1 & \verb|(h_1[1])| & 2 & 1 \\\hline
2 & \verb|(h_2^2[1])| & 7 & 2 \\\hline
3 & \verb|(h_0^3 h_3[1])| & 8 & 4 \\\hline
4 & \verb|(c_0[1])| & 9 & 3 \\\hline
5 & \verb|((Ph_1)[1])| & 10 & 5 \\\hline
\end{tabular}
\caption{\texttt{C2\_AdamsE2\_generators.csv}}\label{tb:gen}
\end{table}

\begin{table}
\begin{tabular}{|c|c|c|}\hline
rel & stem & s\\\hline
0,1,0 & 0 & 1\\\hline
0,1,1;1,2,0 & 2 & 2\\\hline
2,1,1 & 5 & 2\\\hline
0,1,2 & 7 & 3\\\hline
1,1,2;4,1,0 & 8 & 3\\\hline
\end{tabular}
\caption{\texttt{C2\_AdamsE2\_relations.csv}}\label{tb:rel}
\end{table}

\begin{table}
\begin{tabular}{|c|c|c|c|c|}\hline
index & mon & stem & s & d2 \\\hline
0 & 0 & 0 & 0 & \\\hline
0 & 1,1,0 & 1 & 1 & \\\hline
0 & 1 & 2 & 1 & \\\hline
0 & 1,2,0 & 2 & 2 & \\\hline
0 & 2,1,0 & 3 & 1 & \\\hline
\end{tabular}
\caption{\texttt{C2\_AdamsE2\_basis.csv}}\label{tb:basis}
\end{table}

\begin{notation}
    In the table of generators, a row with $id=i$ is the $i$'th generator of the $E_2$-page of \texttt{X}. In this section we will denote this generator by $x_i$ if \texttt{X} is a ring spectrum or $v_i$ if \texttt{X} is a module over some ring spectrum (in our dataset there are only two ring spectra \texttt{S0} and \texttt{tmf}). The columns \texttt{stem,s} are the degrees of the generators. The column \texttt{name} is the latex name we give to the $i$'th generator. However, currently not all generators have latex names. The column \texttt{name} might be sometimes vacant (filled with \verb|[NULL]|).
\end{notation}

\begin{notation}
    In the table of relations, consider the column \verb|rel| for relations. If \verb|X| is a ring spectrum, then \fbox{\texttt{0,1;1,2,5,4}} denotes the relation $$x_0^1+x_1^2x_5^4=0.$$
    If \verb|X| is a module over the ring spectra \verb|R|, then \fbox{\texttt{0,1,3;1,2,5,4,7}} denotes the relation $$x_0^1v_3+x_1^2x_5^4v_7=0$$
    where $v_i$ is the $i$'th generator of \verb|X| and $x_j$ is the $j$'th generator of the ring spectrum \verb|R|.
    In summary, the semicolon can be parsed as the plus sign and the monomials are sequences of sub- and upper-scripts whose lengths are even if \verb|X| is a ring or odd if \verb|X| is a module. The columns \texttt{stem,s} are the degrees of the relations which can be determined by the degrees of the generators.
\end{notation}

\begin{notation}\label{nt:basis}
    In the table of basis, the column \verb|mon| is parsed in the same way as we parse \verb|rel| for relations. The only difference is that each \verb|mon| is a monomial so there should be no semicolon in the column \verb|mon|. The columns \texttt{stem,s} are the degrees of the monomials. The column \texttt{index} is \textbf{important}. When there are multiple monomial basis elements in the same bidegree, they will be indexed by $0,1,\dots,k$. We write
    $$[i]\in \Ext^{\verb|s,stem+s|}(\verb|X|)$$
    to denote the monomial basis element  with index $i$ in $\Ext^{\verb|s,stem+s|}(\verb|X|)$. In the Table of Proofs (see Section \ref{sec:proofs}), the element will be written as
    $$\fbox{\texttt{X (stem, s) [i]}}$$
    in the \verb|info| column. We can also write $[i_1,i_2,...]$ to denote the linear combination of monomials with indices $i_1,i_2,\dots$ and $[]$ to denote the zero element. The \verb|d2| column represents the $d_2$ differentials of the basis elements in the format $[i_1,i_2,...]$ we just described but in degree \verb|(stem-1,s+2)|. If it is empty, it means that the $d_2$ differential is trivial. If it is filled with \verb|[NULL]|, it means that the $d_2$ differential is not calculated yet by the program \texttt{./Adams} using the secondary Steenrod operations.
\end{notation}

Next we explain our naming scheme for the above CW spectra. We chose to use plain text names for these objects because we want to use consistent names for computer codes, command-line arguments for \verb|./Adams| and \verb|./ss|, the website links of plots and the plain text machine proofs. Table \ref{tb:key} is a list of keywords we use in the name and the corresponding elements in $\pi_*(S^0)$ the stable homotopy groups of spheres.

\begin{table}\begin{tabular}{|c|c|}\hline
    keyword & $\pi_*(S^0)$ \\\hline
    \texttt{eta} & $\eta$ \\\hline
    \texttt{nu} & $\nu$\\\hline
    \texttt{sigma} & $\sigma$\\\hline
%    \texttt{2eta} & $2\eta$\\\hline
    \texttt{2nu} & $2\nu$\\\hline
    \texttt{2sigma} & $2\sigma$\\\hline
    % \texttt{etasigma} & $\eta\sigma$\\\hline
    \texttt{sigmasq} & $\sigma^2$\\\hline
    \texttt{theta4} & $\theta_4$\\\hline
    \texttt{theta5} & $\theta_5$\\\hline
    % \texttt{theta5sq} & $\theta_5^2$\\\hline
    % \texttt{etasq} & $\eta^2$\\\hline
    % \texttt{etacube} & $\eta^3$\\\hline
    % \texttt{eta4} & $\eta_4$\\\hline
    % \texttt{eta5} & $\eta_5$\\\hline
    % \texttt{kappa} & $\kappa$\\\hline
    % \texttt{g} & $\bar\kappa$\\\hline
\end{tabular}\caption{Keywords in the names of CW spectra}\label{tb:key}
\end{table}

\begin{notation}
    Let \texttt{a,b,c,...} be keywords, $a,b,c,\dots$ be the corresponding elements in $\pi_*(S^0)$ and $|a|,|b|,|c|$ be the topological degrees of $a,b,c,\dots$. The following examples demonstrate our naming for the CW spectra.
    \begin{itemize}
        \item \verb|S0| is the sphere spectrum $S^0$. \verb|tmf| is the topological modular form $tmf$.
        \item \texttt{Ca} is the cofiber of $a$.
        \item \texttt{CW\_a\_b} is a three-cell complex such that the following two cofiber sequences exist.
        $$\texttt{Ca}\to  \texttt{CW\_a\_b}\to \Sus^{|a|+|b|+2}\verb|S0|$$
        $$\verb|S0|\to  \texttt{CW\_a\_b}\to \Sus^{|a|+1}\texttt{Cb}$$
        The existence of such complex corresponds to $a\cdot b=0$ in $\pi_*(S^0)$.
        \item \texttt{CW\_a\_b\_c} is a four-cell complex such that the following two cofiber sequences exist.
        $$\texttt{CW\_a\_b}\to  \texttt{CW\_a\_b\_c}\to \Sus^{|a|+|b|+|c|+3}\verb|S0|$$
        $$\texttt{Ca}\to  \texttt{CW\_a\_b\_c}\to \Sus^{|a|+|b|+2}\texttt{Cc}$$
        The existence of such complex corresponds to $0\in \langle a, b, c\rangle$ in $\pi_*(S^0)$.
        \item Inductively we can define \verb|CW_a_..._z| as long as it exists.
        \item \verb|Ca_Cb| is the smash product of \texttt{Ca} and \texttt{Cb}. \verb|tmf_X| is the smash product of \verb|tmf| and \verb|X|.
        \item \verb|CW_a_V_b| is a three-cell complex such that the following two cofiber sequences exist.
        $$\verb|Ca|\to \verb|CW_a_V_b|\fto{\text{top cell b}} \Sus^{|b|+1}\verb|S0|$$
        $$\verb|Cb|\to \verb|CW_a_V_b|\fto{\text{top cell a}} \Sus^{|a|+1}\verb|S0|$$
        \item \verb|DX| is a \emph{$k$-fold suspension} of the Spanier–Whitehead dual of \verb|X| such that the lowest dimension of the cells is $0$.
        \item \verb|CW_a_A_b| is equivalent to \verb|DCW_a_V_b|.
        \item \verb|C2h4| is the same as \verb|CW_2_sigmasq|. \verb|C2h5| is the same as \verb|CW_2_theta4|.  \verb|C2h6| is the same as \verb|CW_2_theta5|.
        \item \verb|CW_a_b_V_c| is a four-cell complex such that the following two cofiber sequences exist.
        $$\verb|CW_a_b|\to \verb|CW_a_b_V_c|\to \Sus^{|c|+1}\verb|S0|$$
        $$\verb|Cc|\to \verb|CW_a_b_V_c|\to \Sus^{|a|+1}\verb|Cb|$$
        \item \verb|CW_a_b_c_Eq_d_e| is a five-cell complex such that the following two cofiber sequences exist.
        $$\verb|CW_a_b_V_d|\to \verb|CW_a_b_c_Eq_d_e|\to \Sus^{|a|+|b|+|c|+3}\verb|S0|$$
        $$\verb|Ca|\to \verb|CW_a_b_c_Eq_d_e|\to \Sus^{\min\{|a|+|b|+2,|d|+1\}}\verb|CW_c_A_e|$$
        The existence of such complex corresponds to $d\cdot e\in \langle a, b, c\rangle$ in $\pi_*(S^0)$.
        \item \verb|Joker| is equivalent to \verb|CW_2_eta_2_Eq_eta_eta|.
        \item \verb|RPm_n| is equivalent to the stunted real projective space $\bR P_m^n$ for positive integers \verb|m,n|.
        % \item \verb|CPm_n| is equivalent to the stunted complex projective space $\bC P_m^n$ for positive integers \verb|m,n|.
        \item \verb|Fphi| is the suspension of $F\phi$ the fiber of the Kahn-Priddy map $\phi: \bR P_1^\infty\to S^0$ (or the cofiber of $\phi$).
        \item \verb|Fphik| is the suspension of the $k$-skeleton of $F\phi$ for an integers \verb|k|.
    \end{itemize}
\end{notation}
\begin{remark}
    All the CW spectra we considered are connective spectra. The current programs only work well with connective spectra.
\end{remark}
\begin{remark}
    If there are two or more CW spectra with different homotopy types that satisfy the same naming condition, then the name represents all choices and we do not assume which one we pick. The current important input of the program is the existence of these cofiber sequences and the ones derived from these using the properties of triangulated categories. We will make assumptions on homotopy classes in the future when a new feature of the program is complete and it will enable the \verb|./ss| program to rigorously determine if there exists at least one homotopy class that satisfies multiple constraints.
\end{remark}

All the CW spectra we consider are known to exist in the literature. One major example is the existence of $tmf$ and its 15-skeleton 
$$\verb|CW_sigma_nu_eta_2|,$$
and we can get many subquotients of $\verb|CW_sigma_nu_eta_2|$ and duals of these such as $\verb|CW_2_eta_nu|$. The existence of $\verb|C2h4|, \verb|C2h5|, \verb|C2h6|$ is due to the strong Kervaire invariant problem for $\sigma^2,\theta_4,\theta_5$ (See \cite{BarrattJonesMahowald} \cite{Xu}).

\subsection{Maps between Adams $E_2$-pages}
Tables \ref{tb:map1}--\ref{tb:map2} give the list of maps between CW spectra we have considered and the maximum computed internal degrees of these maps. The maximum computed internal degree of a map $f: X\to \Sus^{k}Y$ refers to the maximum $T$ such that such that the induced map $$f:E_2^{s,t}(X)\to E_2^{s+\AF(f),t-k+\AF(f)}(Y)$$ is computed in the range $t\le T$.\vspace{2pt}

\input{cw_map_data}

For a map \texttt{X\_\_Y} (we will explain the naming later), the data can be found in the database file \texttt{map\_AdamsSS\_X\_to\_Y.db} in \texttt{kervaire\_database.zip} or in the plain-text table file \texttt{map\_X\_to\_Y.csv}. (A few filenames created by a newer version of \verb|./Adams| drop `\verb|to|' from the names.) Table \ref{tb:mapcsv} lists some rows from \verb|map_C2_to_C2h4.csv| to illustrate the following notation.

\begin{notation}
    In the table of maps, there are only two columns \verb|id,map|. The \verb|map| column gives us the image of the \verb|id|'th generator from $E_2^{*,*}(X)$ to $E_2^{*,*}(Y)$, which is parsed in the same way as we parse \verb|rel| in the table of relations. For example, we can read from Table \ref{tb:mapcsv} that the image of $v_{36}\in E_2^{*,*}(\verb|C2|)$ is
    $$v_{63}+x_1x_{48}v_0\in E_2^{*,*}(\verb|C2h4|).$$
\end{notation}

\begin{table}
\begin{tabular}{|c|c|}\hline
id & map \\\hline
34 & 61 \\\hline
35 & 62 \\\hline
36 & 63;1,1,48,1,0 \\\hline
37 & 66;2,1,53 \\\hline
38 & 67 \\\hline
\end{tabular}
\caption{\texttt{map\_C2\_to\_C2h4.csv}}\label{tb:mapcsv}
\end{table}

We explain our naming scheme for these maps.
\begin{notation}
    Let \verb|X,Y,Z| be names of CW spectra. The following examples demonstrate our naming for the maps between CW spectra.
    \begin{itemize}
        \item If \verb|X| has strictly fewer cells than \verb|Y| and there exists a unique $k$ and a unique inclusion (up to homotopy) such that \verb|X| is a subcomplex of $\Sus^k$\verb|Y|, then \verb|X__Y| denotes this subcomplex map \verb|X|$\hookrightarrow\Sus^k$\verb|Y|.
        \item If \verb|X| has strictly more cells than \verb|Y| and there exists a unique $k$ and a unique quotient such that $\Sus^k$\verb|Y| is a quotient complex of \verb|X|, then \verb|X__Y| denotes this quotient map \verb|X|$\to\Sus^k$\verb|Y|.
        \item If we have a cofiber sequence
        $$\verb|X|\to \Sus^m\verb|Y|\to \Sus^n\verb|Z|\to \Sus\verb|X|$$
        whose first map is a subcomplex inclusion denoted by \verb|X__Y| and whose second map is a quotient map denoted by \verb|Y__Z|, then the third boundary map is denoted by \verb|Z__Q_Y|.
        \item For key words \verb|a,b,c,d,e|, the map \verb|CW_d_e__CW_a_b_by_c| is equivalent to 
        $$\verb|CW_d_e__Q_CW_a_b_c_d_e|.$$
        Similarly, if we have a map with name \verb|X__Y_by_c| it means that \verb|X| has a unique bottom cell and \verb|Y| has a unique top cell and the map attaches the bottom cell of \verb|X| to the top cell of \verb|Y| by \verb|c|. In other words, these two cells form a subquotient of the cofiber of the map which is equivalent to \verb|Cc|.
        \item For numbers \verb|n,l,p,q| where $p<q=n-1<n<l$, the map \verb|RPn_l__RPp_q| is equivalent to \verb|RPn_l__Q_RPp_l|. In other words, the following three maps form a cofiber sequence.
        $$(\verb|RPp_q__RPp_l|,~\verb|RPp_l__RPn_l|,~\verb|RPn_l__RPp_q|)$$
        \item In plain-text names, we use the letter \verb|m| to denote the negative sign. For example, \verb|RP1_256__RPm7_0| corresponds to the map
        $$\bR P_1^{256}\to \Sus \bR P_{-7}^0\simeq \Sus^{-7}\bR P_1^8.$$
        \item Note that $\bR P_0^n\equiv S^0\vee \bR P_1^n$. We use \verb|RP1_256__S0| to denote the truncated Kahn-Priddy map, which is a composition of the following maps.
        $$\bR P_1^{256}\to \bR P_0^{256}\to \Sus \bR P_{-1}^{-1}=S^0.$$
        \item \verb|S0__tmf| and \verb|X__tmf_X| are the Hurewicz maps.
    \end{itemize}
\end{notation}

All the maps except the Hurewicz maps \verb|S0__tmf| and \verb|X__tmf_X| we considered are maps in the cofiber sequences generated by inclusions $X\subset Y$ of CW complexes, i.e., the homology maps they induce are injective, surjective or trivial.

\subsection{Adams $d_2$ differentials}
Table \ref{tb:d2} demonstrates the range of $d_2$ differentials we have computed and used in the computer programs. The maximum internal degree refers to the maximum $T$ such that all the $d_2$ differentials of $x\in E_2^{s,t\le T}$ are known for the corresponding CW spectra.
\begin{table}
\begin{tabular}{|l|c|}\hline
$X$ & max(t)\\\hline
\verb|S0| & 176 \\\hline
\verb|C2| & 198 \\\hline
\verb|Ceta| & 169 \\\hline
\verb|Cnu| & 176 \\\hline
\verb|Csigma| & 180 \\\hline
\verb|C2h4| & 167 \\\hline
\verb|C2h5| & 170 \\\hline
\verb|C2h6| & 176 \\\hline
\verb|DC2h4| & 134 \\\hline
\verb|DC2h5| & 138 \\\hline
\verb|DC2h6| & 145 \\\hline
\verb|CW_2_eta| & 159 \\\hline
\verb|CW_2_eta_nu| & 149 \\\hline
\verb|CW_2_eta_nu_sigma| & 149 \\\hline
% \verb|CW_2_eta_theta5sq| & 156 \\\hline
\verb|CW_eta_2| & 144 \\\hline
\verb|CW_eta_nu| & 149 \\\hline
\verb|CW_eta_nu_sigma| & 149 \\\hline
\verb|CW_nu_eta| & 149 \\\hline
\verb|CW_nu_eta_2| & 149 \\\hline
\verb|CW_nu_sigma| & 149 \\\hline
\verb|CW_sigma_nu| & 148 \\\hline
\verb|CW_sigma_nu_eta| & 149 \\\hline
\verb|CW_sigma_nu_eta_2| & 149 \\\hline
\verb|Joker| & 135 \\\hline
\verb|RP1_4| & 135 \\\hline
\verb|RP1_6| & 128 \\\hline
\verb|RP1_8| & 125 \\\hline
\end{tabular}\caption{Computed $d_2$ differentials}\label{tb:d2}
\end{table}\vspace{5pt}

We refer the reader to \cite{Chua} and \cite{Nassau} for the algorithm that computes $d_2$ differentials using the secondary Steenrod operations encoded in a differential graded algebra. The first author also implement this algorithm in the program \verb|./Adams|. The $d_2$ data we computed is already contained the table of basis (see Notation \ref{nt:basis}).

\section{Adams spectral sequences and Extension spectral sequences}\label{sec:diff}
In this section we describe how the computer program \verb|./ss| organizes the data of the Adams spectral sequence and the extension spectral sequence and how to read the differentials and extensions generated by the computer program.

The program \texttt{./ss} organizes the Adams spectral sequence by the following.
$$0=B_1^{s,t}\subseteq B_2^{s,t}\subseteq B_3^{s,t}\subseteq \cdots\subseteq Z_3^{s,t}\subseteq Z_2^{s,t}\subseteq E_2^{s,t}$$
where $B_r$ is the subgroup of the image of $d_r$ differentials and $Z_r$ is the subgroup of $d_r$-cycles. We can interpret Adams differentials as isomorphisms
\begin{equation}\label{eq:dr}
    d_r: Z_{r-1}^{s,t}/Z_r^{s,t}\fto{~~\iso~~} B_r^{s+r,t+r-1}/B_{r-1}^{s+r,t+r-1}
\end{equation}
and therefore we can also consider the inverse
\begin{equation}\label{eq:drinv}
    d_r^{-1}: B_{r}^{s,t}/B_{r-1}^{s,t}\fto{~~\iso~~} Z_{r-1}^{s-r,t-r+1}/Z_r^{s-r,t-r+1}.
\end{equation}
Table \ref{tb:dabf1541} is an illustration of how we store the spectral sequence in memory. In Table \ref{tb:dabf1541} we can directly see that $B_2^{11,123+11}(S^0)=B_3=B_4$ is generated by $h_0^2x_{123,9}$, $B_5^{11,123+11}(S^0)=B_6$ is generated by $h_0^2x_{123,9}$ and $h_5x_{92,10}$, and $Z_6^{11,123+11}(S^0)=Z_5=Z_4=Z_3$ is generated by $h_0^2x_{123,9}$, $h_5x_{92,10}$ and $x_{123,11,2}+x_{123,11}+h_0h_6[B_4]$, etc.

\begin{table}\scalebox{0.85}{\begin{tabular}{|c|l|l|l|}\hline
    \multirow{5}{*}{11} & $h_0^2x_{123,9}$ & $d_{2}^{-1}$ & $x_{124,9}$ \\\cline{2-4}
     & $h_5x_{92,10}$ & $d_{5}^{-1}$ & $x_{124,6}$ \\\cline{2-4}
     & $x_{123,11,2}+x_{123,11}+h_0h_6[B_4]$ & $d_{7}$ & $h_1x_{121,17}$ \\\cline{2-4}
     & $x_{123,11}+h_0h_6[B_4]$ & $d_{3}$ & $h_0x_{122,13}$ \\\cline{2-4}
     & $h_0h_6[B_4]$ & $d_{2}$ & $h_0^2h_6Md_0$ \\\hline\hline
    \multirow{3}{*}{10} & $h_0x_{123,9}$ & $d_{2}^{-1}$ & $x_{124,8}$ \\\cline{2-4}
     & $x_{123,10}+h_6[B_4]$ & $d_{3}^{-1}$ & $x_{124,7}$ \\\cline{2-4}
     & $h_6[B_4]$ & $d_{2}$ & $h_0h_6Md_0$ \\\hline\hline
    \multirow{2}{*}{9} & $x_{123,9}+h_0x_{123,8}$ & $d_{12}$ & $?$ \\\cline{2-4}
     & $h_0x_{123,8}$ & $d_{3}$ & $h_0x_{122,11}+h_0h_6Md_0$ \\\hline\hline
    \multirow{1}{*}{8} & $x_{123,8}$ & $d_{3}$ & $x_{122,11}+h_6Md_0$ \\\hline
\end{tabular}}\caption{$E_2^{*,*}(S^0)$ \text{stem}=123}\label{tb:dabf1541}
\end{table}

We write differentials in terms of elements in $E_2^{*,*}$. The value of $d_r^{\pm 1}$ should be viewed as a coset. By (\ref{eq:dr}) and (\ref{eq:drinv}) we see that
\begin{itemize}
    \item if we get $d_r(x)=y_1$ and $d_r(x)=y_2$ for two different reasons, then we have $y_1-y_2\in B_{r-1}$;
    \item if we get $d_r^{-1}(x)=z_1$ and $d_r^{-1}(x)=z_2$ for two different reasons, then we have $z_1-z_2\in Z_r$.
\end{itemize}

Consider a cofiber sequence
$$X\fto{f} Y\fto{g} Z.$$
The program \texttt{./ss} organizes the extension spectral sequence (see \cite{LWX126}) by the following.
$${B}_\infty^{s,t}\subseteq \BESS{f}_0^{s,t}\subseteq \BESS{f}_1^{s,t}\subseteq \cdots \subseteq \ZESS{g}_1^{s,t} \subseteq \ZESS{g}_0^{s,t} \subseteq {Z}_\infty^{s,t}.$$
where $B_r$ is the subgroup of the image of $d_r$ differentials and $Z_r$ is the subgroup of $d_r$-cycles. We can interpret extensions as isomorphisms
\begin{equation}\label{eq:drf}
    d_r^f: \ZESS{f}_{r-1}^{s,t}(X)/\ZESS{f}_r^{s,t}(X)\fto{~~\iso~~} \BESS{f}_r^{s+r,t+r}(Y)/\BESS{f}_{r-1}^{s+r,t+r}(Y)
\end{equation}
when $r>0$ and
\begin{equation}\label{eq:drfinv}
    d_0^f: Z_\infty^{s,t}(X)/\ZESS{f}_0^{s,t}(X)\fto{~~\iso~~} \BESS{f}_0^{s+r,t+r}(Y)/B_{\infty}^{s+r,t+r}(Y)
\end{equation}
when $r=0$.

We write differentials in terms of elements in $E_2^{*,*}$. The value of $(d_r^f)^{\pm 1}$ should be viewed as a coset. By (\ref{eq:drf}) and (\ref{eq:drfinv}) we see that
\begin{itemize}
    \item if $r>0$ and we get $d_r^f(x)=y_1$ and $d_r^f(x)=y_2$ for two different reasons, then we have $y_1-y_2\in \BESS{f}_{r-1}(Y)$;
    \item if $r=0$ and we get $d_0^f(x)=y_1$ and $d_0^f(x)=y_2$ for two different reasons, then we have $y_1-y_2\in B_\infty(Y)$.
    \item we leave it to the reader to think about what happens to $(d^f_r)^{-1}$ in similar situations.
\end{itemize}

The program \verb|./ss| can compute a lot differentials and extensions based on the data structure we describe above and various theorems such as the Leibniz rule and naturality. The proofs are given in the Table of Proofs (Section \ref{sec:proofs}). There are three manually added differentials: $d_5 h_0^{24}h_6=h_0^2P^6d_0$ and
$d_6 h_0^{55}h_7=h_0^2x_{126,60}$ in the Adams spectral sequence of $S^0$ (from the image of $J$), and $d_3 v_2^{16}=\beta^5g$ in the Adams spectral sequence of $\mathit{tmf}$, derived from power operations (Bruner--Rognes \cite{BR21}).

The Adams differentials for a spectrum \verb|X| produced by the machine are available in the database \verb|X_AdamsSS.db| within \verb|kervaire_database.zip|. Alternatively, the data can be accessed in the plain-text table file \verb|X_AdamsE2_ss.csv| within \verb|kervaire_csv.rar|.

The extensions in a cofiber sequence \verb|X__Y__Z| are available in the database \verb|cofseq_X__Y__Z.db| within \verb|kervaire_database.zip|. Alternatively, the data can be accessed in the plain-text table file \verb|cofseq_X__Y__Z.csv| within \verb|kervaire_csv.rar|.

\begin{table}
\begin{tabular}{|l|}\hline
\verb|C2__C2_C2__C2|\\\hline
\verb|C2__C2_Ceta__C2|\\\hline
\verb|C2__C2h4__S0|\\\hline
\verb|C2__C2h5__S0|\\\hline
\verb|C2__C2h6__S0|\\\hline
\verb|C2__CW_2_V_eta__S2|\\\hline
\verb|C2__CW_2_eta__S0|\\\hline
\verb|C2__CW_2_eta_nu__Cnu|\\\hline
\verb|C2__CW_2_eta_nu_sigma__CW_nu_sigma|\\\hline
\verb|CW_2_eta__CW_2_eta_nu__S0|\\\hline
\verb|CW_2_eta__CW_2_eta_nu_sigma__Csigma|\\\hline
\verb|CW_2_eta__Joker__Ceta|\\\hline
\verb|CW_2_eta_nu__CW_2_eta_nu_sigma__S0|\\\hline
\verb|CW_eta_2__RP1_6__CW_2_eta|\\\hline
\verb|CW_eta_nu__CW_eta_nu_sigma__S0|\\\hline
\verb|CW_nu_eta__CW_nu_eta_2__S0|\\\hline
\verb|Ceta__C2_Ceta__Ceta|\\\hline
\verb|Ceta__CW_2_V_eta__S1|\\\hline
\verb|Ceta__CW_eta_2__S0|\\\hline
\verb|Ceta__CW_eta_nu__S0|\\\hline
\verb|Ceta__CW_eta_nu_sigma__Csigma|\\\hline
\verb|Ceta__Joker__CW_eta_2|\\\hline
\verb|Cnu__CW_nu_eta_2__C2|\\\hline
\verb|Cnu__CW_nu_eta__S0|\\\hline
\verb|Cnu__CW_nu_sigma__S0|\\\hline
\verb|Csigma__CW_sigma_2sigma__S0|\\\hline
\verb|Csigma__CW_sigma_nu__S0|\\\hline
\verb|Csigmasq__DC2h4__S0|\\\hline
\verb|Ctheta4__DC2h5__S0|\\\hline
\verb|Ctheta5__DC2h6__S0|\\\hline
\verb|Fphi__RP1_256__S0|\\\hline
\verb|RP1_2__RP1_256__RP3_256|\\\hline
\verb|RP1_2__RP1_4__RP3_4|\\\hline
\verb|RP1_2__RP1_6__RP3_6|\\\hline
\verb|RP3_4__RP3_6__RP5_6|\\\hline
\verb|S0__C2__S0|\\\hline
\verb|S0__C2h4__Csigmasq|\\\hline
\verb|S0__C2h5__Ctheta4|\\\hline
\verb|S0__C2h6__Ctheta5|\\\hline
\verb|S0__C2sigma__S0|\\\hline
\verb|S0__CW_2_A_eta__C2|\\\hline
\verb|S0__CW_2_eta__Ceta|\\\hline
\verb|S0__CW_2_eta_nu__CW_eta_nu|\\\hline
\verb|S0__CW_2_eta_nu_sigma__CW_eta_nu_sigma|\\\hline
\end{tabular}\caption{Cofiber sequences in the database}\label{tb1:cofseq}
\end{table}\vspace{5pt}
\begin{table}
\begin{tabular}{|l|}\hline
\verb|S0__CW_eta_2__C2|\\\hline
\verb|S0__CW_eta_nu__Cnu|\\\hline
\verb|S0__CW_eta_nu_sigma__CW_nu_sigma|\\\hline
\verb|S0__CW_nu_eta_2__CW_eta_2|\\\hline
\verb|S0__CW_nu_eta__Ceta|\\\hline
\verb|S0__CW_nu_sigma__Csigma|\\\hline
\verb|S0__CW_sigma_nu__Cnu|\\\hline
\verb|S0__Ceta__S0|\\\hline
\verb|S0__Cnu__S0|\\\hline
\verb|S0__Csigma__S0|\\\hline
\verb|S0__Csigmasq__S0|\\\hline
\verb|S0__Ctheta4__S0|\\\hline
\verb|S0__Ctheta5__S0|\\\hline
\verb|S0__DC2h4__C2|\\\hline
\verb|S0__DC2h5__C2|\\\hline
\verb|S0__DC2h6__C2|\\\hline
\verb|S1__CW_2_A_eta__Ceta|\\\hline
\end{tabular}\caption{Cofiber sequences in the database}\label{tb2:cofseq}
\end{table}\vspace{5pt}

For illustration, Table~\ref{tb:ss} is a table of some rows from the table of Adams spectral sequence \texttt{S0\_AdamsE2\_ss.csv} while Table~\ref{tb:cs} is a table of some rows from the table of extension spectral sequence \texttt{cofseq\_S0\_\_C2\_\_S0.csv}.

\begin{table}
\begin{tabular}{|c|c|c|c|c|}\hline
stem & s & base & diff & level \\\hline
15 & 1 & 0 & 0 & 9998 \\\hline
0 & 17 & 0 & \verb|[NULL]| & 9000 \\\hline
11 & 6 & 0 & \verb|[NULL]| & 9000 \\\hline
14 & 3 & 0 & 0 & 2 \\\hline
15 & 2 & 0 & 0 & 9997 \\\hline
\end{tabular}
\caption{\texttt{S0\_AdamsE2\_ss.csv}}\label{tb:ss}
\end{table}

\begin{table}
\begin{tabular}{|c|c|c|c|c|c|}\hline
iC & stem & s & base & diff & level\\\hline
1 & 87 & 42 & 0 & 1 & 1\\\hline
1 & 127 & 3 & 0 & \verb|[NULL]| & 9996\\\hline
1 & 126 & 7 & 1 & \verb|[NULL]| & 9998\\\hline
1 & 118 & 19 & 0 & 1,3 & 10000\\\hline
1 & 134 & 4 & 0 & \verb|[NULL]| & 9996\\\hline
1 & 118 & 21 & 0 & 0 & 10000\\\hline
\end{tabular}
\caption{\texttt{cofseq\_S0\_\_C2\_\_S0.csv}}\label{tb:cs}
\end{table}

\begin{notation}\label{nt:4b392f50}
    In the table of Adams spectral sequence, columns \verb|stem,s| refers to the degrees of the elements in the \verb|base| column. The columns \verb|base,diff| should be parsed as linear combinations of the monomial basis in the corresponding bidegree. If \verb|level=r<5000|, the \verb|base| column is hit by a $d_r$-differential from \verb|diff|. If \verb|level=10000-r>5000|, the \verb|base| column supports a $d_r$-differential hitting \verb|diff|. If \verb|level=9000|, it means that \verb|base| is a permanent cycle. If \verb|diff| is filled with \verb|[NULL]|, it means that the value of \verb|diff| is undetermined. The reader can see that we assume that all visible nontrivial differentials in the data are of length shorter than 1000.
\end{notation}

\begin{notation}
    In the table of extension spectral sequence corresponding to a cofiber sequence
    $$X\fto{f} Y\fto{g} Z,$$
    the column \verb|iC=0,1,2| refers to \verb|X,Y,Z| where \verb|base| lives respectively. The other columns are similar to those in Notation \ref{nt:4b392f50}. If \verb|iC=1| and \verb|level=r<5000|, then \verb|base| is hit by an $f$-extension which jumps filtration by $r$. If \verb|iC=1| and \verb|level=10000-r>5000|, then \verb|base| supports an $g$-extension which jumps filtration by $r$.
\end{notation}

\begin{remark}
    All the Adams differentials and extensions can be also visually accessed from the website \cite{LinPlot}.
\end{remark}

\section{Proofs}\label{sec:proofs}
In this section we explain the structure and notations in our machine generated proofs. We call it the Table of Proofs. The Table of Proofs is provided in two different forms in the Zenodo repo: a single sqlite3 database file \verb|proofs.db| and 22 plain-text tables \verb|proofs-part1.csv|--\verb|proofs-part22.csv|. The plain-text table is divided into 22 table files because there are more than 21 million rows in the table and some table viewers like Microsoft Excel has a row limit of one million. We encourage the reader to use the open source software DB Browser for SQLite to access the table in \verb|proofs.db|. This tool can help you quickly search for what you want in big tables. You can also execute advanced SQL queries in the software.

The following are two typical rows in the Table of Proofs.
\begin{table}[h]\scalebox{0.85}{\begin{tabular}{|l|l|l|l|l|l|l|l|l|l|l|}\hline
id & depth & reason & name & stem & s & t & r & x & dx & info\\\hline
248925 & 0 & D & \texttt{Ceta} & 107 & 16 & 123 & 3 & 0,2 & 0,1 & \texttt{NULL}\\\hline
2042139 & 0 & D & \texttt{Cnu\_\_CW\_nu\_eta\_2\_\_C2:1} & 15 & 2 & 17 & 2 & 0 & 0 & \texttt{NULL}\\\hline
\end{tabular}}\caption{Two rows from the Table of Proofs}\label{tb:proofs}
\end{table}\vspace{5pt}

The first column \verb|id| refers to the line number in the table. The columns \verb|depth| and \verb|reason| refer to the structure of the proof which we will explain in more details later. The column \verb|name| refers to the name of a CW spectra or a map in a cofiber sequence. 
\begin{itemize}
    \item If \verb|name| is the name of a CW spectrum like the first row, then the columns \verb|stem,s,r,x,dx| refer to an Adams differential of the spectrum. In our example, this is $$d_3([0,2])=[0,1]$$ for $$[0,2]\in E_2^{16,107+16}(\verb|Ceta|)$$ where $[0,2]$ denotes the linear combination of the 0th and 2nd basis elements in $E_2^{16,107+16}$ in the basis table of \verb|Ceta| and $[0,1]$ denotes the linear combination in the bidegree that this differential hits.
    \item If \verb|name| is the name of a map in a cofiber sequence, it is of the form \verb|X__Y__Z:n| where
    $$\verb|X|\fto{f} \Sus^?\verb|Y|\fto{g} \Sus^? \verb|Z|\fto{h} \Sus \verb|X|$$
    is a cofiber sequence and $n=0,1,2$ correspond to $f,g,h$ respectively. The columns \verb|stem,s,r,x,dx| refer to an extension of the map. In our example (second row of Table \ref{tb:proofs}), consider the cofiber sequence
    $$\verb|Cnu|\fto{f} \verb|CW_nu_eta_2|\fto{g} \Sus^{6}\verb|C2|\fto{h} \Sus \verb|Cnu|.$$
    The second row of the table indicates that there is a $g$-extension from the permanent cycle $$[0]\in E_2^{2,15+2}(\verb|CW_nu_eta_2|)$$ to the permanent cycle $$[0]\in E_2^{4,15+4-6}(\verb|C2|)$$
    which jumps filtration by $r=2$.
\end{itemize}

Next we explain the structure of the Table of Proofs. A typical proof of an Adams differential has the following form.

\begin{proofstr}\label{ps:e57cb438}
Assume that the value of an Adams differential $d_r(x)$ has four candidates $y_1,y_2,y_3,y_4$. (Usually the number of candidates is a power of two.)
\begin{itemize}
    \item If $d_r(x)=y_1$, propagate (by the Leibniz rule, naturality, ...). Get a contradiction.
    \item If $d_r(x)=y_2$, propagate. Get a contradiction.
    \item If $d_r(x)=y_3$, propagate. Get a contradiction.
\end{itemize}
Then we can conclude that $d_r(x)=y_4$.
\end{proofstr}

The following 4 rows from the Table of Proofs is a concrete form of Proof Structure \ref{ps:e57cb438}.\vspace{5pt}
\begin{center}\scalebox{1.0}{\begin{tabular}{|l|l|l|l|l|l|l|l|l|l|l|}\hline
id & depth & reason & name & stem & s & t & r & x & dx & info\\\hline
325477 & 1 & T & Csigma & 116 & 10 & 126 & 4 & 1 & 3 & 1,2  \texttt{proof1}\\\hline
325478 & 1 & T & Csigma & 116 & 10 & 126 & 4 & 1 & 0 & 0,2  \texttt{proof2}\\\hline
325479 & 1 & T & Csigma & 116 & 10 & 126 & 4 & 1 & 0,3 & 0,1  \texttt{proof3}\\\hline
325480 & 0 & D & Csigma & 116 & 10 & 126 & 4 & 1 &  & \texttt{NULL}\\\hline
\end{tabular}}\end{center}\vspace{5pt}

Here the \verb|depth| column can be thought as the nesting level of the latex \verb|itemize| environment. Rows with \verb|reason=`T'| correspond to assumptions (`T' is a symbol for `Try'). Each assumption row has a multi-line proof in the \verb|info| column which explains why the assumption leads to a contradiction. For example, the actual content of \verb|proof1| is the following.\vspace{5pt}
\begin{mdframed}\begin{verbatim}
Get Csigma (116,10) d_3[1]=[1,2]. Apply the Leibniz rule with 
S0 (0,1) d_3[0]=[] and get Csigma (116,11) d_3[]=[1,2].
However, Csigma (115,14) [1,2] is not in B_2.
\end{verbatim}\end{mdframed}\vspace{5pt}
In \verb|proof1| the element \verb|Csigma (115,14) [1,2]| refers to the element
$$[1,2]\in E_2^{14,115+14}(\verb|Csigma|).$$
Rows with \verb|reason=`D'| correspond to the conclusions after we rule out enough candidates (`T' rows) of the value of Adams differentials or extensions (`D' is a symbol for `Deduction').

Microsoft Excel tip: Since the text in the \verb|info| column usually contains multiple lines, users who view the \verb|csv| files with Excel should adjust the height of the edit box indicated by the following screenshot in order to view the whole content of the text.
\begin{center}
\includegraphics[width=12cm]{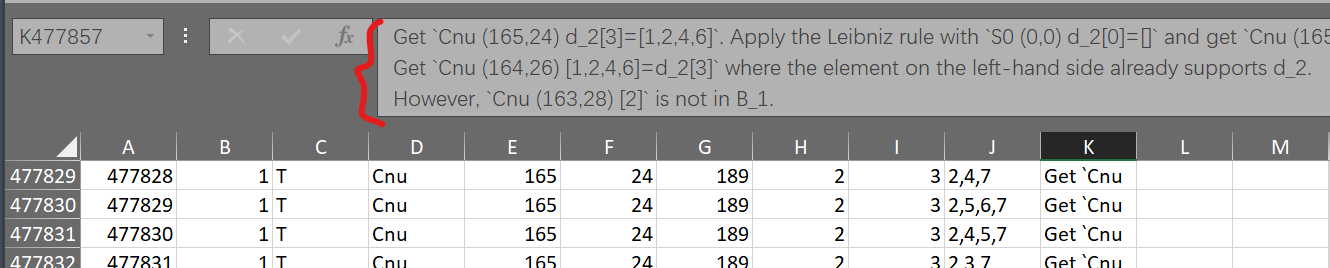}
\end{center}\vspace{5pt}

Next we give examples of more complex proofs of Adams differentials.
\begin{proofstr}\label{ps:cfb2f966}
Assume that the value of an Adams differential $d_r(x)$ has four candidates $0, y_1,y_2,y_3$.
\begin{itemize}
    \item If $d_r(x)=0$, propagate. Assume that $d_{r+1}(x)$ has four possible values $z_1, z_2, z_3, z_4$.
    \begin{itemize}
        \item Assume $d_{r+1}(x)=z_1$. Propagate. Get a contradiction.
        \item Assume $d_{r+1}(x)=z_2$. Propagate. Get a contradiction.
        \item Assume $d_{r+1}(x)=z_3$. Propagate. Get a contradiction.
        \item Assume $d_{r+1}(x)=z_4$. Propagate. Get a contradiction.
    \end{itemize}
    \item If $d_r(x)=y_1$, propagate. Get a contradiction.
    \item If $d_r(x)=y_2$, propagate. Get a contradiction.
\end{itemize}
Then we can conclude that $d_r(x)=y_3$.
\end{proofstr}

The following 8 rows from the Table of Proofs is an example of Proof Structure \ref{ps:cfb2f966}.\vspace{5pt}
\begin{center}\scalebox{1.0}{\begin{tabular}{|l|l|l|l|l|l|l|l|l|l|l|}\hline
id & depth & reason & name & stem & s & t & r & x & dx & info\\\hline
2463208 & 1 & T & \texttt{CW\_nu\_sigma} & 112 & 12 & 124 & 4 & 0 &  & \texttt{proof1}\\\hline
2463209 & 2 & T & \texttt{CW\_nu\_sigma} & 112 & 12 & 124 & 5 & 0 &  & \texttt{proof2}\\\hline
2463210 & 2 & T & \texttt{CW\_nu\_sigma} & 112 & 12 & 124 & 5 & 0 & 3 & \texttt{proof3}\\\hline
2463211 & 2 & T & \texttt{CW\_nu\_sigma} & 112 & 12 & 124 & 5 & 0 & 0 & \texttt{proof4}\\\hline
2463212 & 2 & T & \texttt{CW\_nu\_sigma} & 112 & 12 & 124 & 5 & 0 & 0,3 & \texttt{proof5}\\\hline
2463213 & 1 & T & \texttt{CW\_nu\_sigma} & 112 & 12 & 124 & 4 & 0 & 0 & \texttt{proof6}\\\hline
2463214 & 1 & T & \texttt{CW\_nu\_sigma} & 112 & 12 & 124 & 4 & 0 & 0,3 & \texttt{proof7}\\\hline
2463215 & 0 & D & \texttt{CW\_nu\_sigma} & 112 & 12 & 124 & 4 & 0 & 3 & \texttt{NULL}\\\hline
\end{tabular}}\end{center}\vspace{5pt}\vspace{5pt}

The reader can imagine that there could be proofs with even more complex structure. Consider the following.

\begin{proofstr}\label{ps:833b7e80}
Assume that the value of an Adams differential $d_r(x)$ has four candidates $0, y_1,y_2,y_3$.
\begin{itemize}
    \item If $d_r(x)=0$, propagate. Assume that $d_{r+1}(x)$ has two possible values $z,0$.
    \begin{itemize}
        \item Assume $d_{r+1}(x)=z$. Propagate. Get a contradiction.
    \end{itemize}
    We get $d_{r+1}(x)=0$. Assume that $d_{r+2}(x)$ has two possible values $w_1,w_2$.
    \begin{itemize}
        \item Assume $d_{r+2}(x)=w_1$. Propagate. Get a contradiction.
        \item Assume $d_{r+2}(x)=w_2$. Propagate. Get a contradiction.
    \end{itemize}
    \item If $d_r(x)=y_1$, propagate. Get a contradiction.
    \item If $d_r(x)=y_2$, propagate. Get a contradiction.
\end{itemize}
Then we can conclude that $d_r(x)=y_3$.
\end{proofstr}

A proof with the above structure would be 8 rows from the Table of Proofs with (depth, reason)=(1,T), (2,T), (1,D), (2,T), (2,T), (1,T), (1,T), (0,D). If at some point in the proof we got a permanent cycle, we can even make assumptions about the value of extensions. The reader can find such examples in the Table of Proofs.

Next we describe other types of rows that can be computed or derived from former differentials labeled with different \verb|reason|s.
\begin{itemize}
    \item Rows with \verb|reason=`d2'|. This means that the $d_2$ differential is calculated by \verb|./Adams| using secondary Steenrod operation algorithm.
    \item Rows with \verb|reason=`N'|. This means that we get a differential from a previous one by naturality.
    \item Rows with \verb|reason=`G'|. This means that we get a trivial differential for degree reason.
    \item Rows with \verb|reason=`XX'|. This corresponds to the case that if we have $d_{r-1}(x)=0$ then we automatically get $d_r(x^2)=0$.
    \item Rows with \verb|reason=`XY'|. This corresponds to the case that if we have $d_{r-1}(x)=0$ then no matter what $d_r(x)$ is (when there are two candidates or more), we always get the same value for $d_r(xy)$ by Leibniz rule or $d_r(f(x))$ by naturality.
    \item Rows with \verb|reason=`ToCs'|. This means that in the extension spectral sequence, we have tried all possible values of an extension (linear combinations of elements in the $E_2$-page) and the outcome is that certain element in the Adams spectral sequence should survive to the $E_\infty$-page and it becomes a summand of the value of the extension.
    \item Rows with \verb|reason=`OutCsI'|. This means that in the extensions spectral sequence, two permanent cycles should be the same on the $E_\infty$-page. Hence the difference \verb|dx| should be hit by an Adams differential. In those rows the columns \verb|stem,s,t| refer to the degree of \verb|dx| in stead of \verb|x|.
    \item Rows with \verb|reason=`CsCm'|. These are extensions obtained from some simple homotopy relations. For example, if $x\in E_\infty(S^0)$ is the target of an $\eta$-extension, then it should detect a $\nu$-torsion ($d_n^{\nu}(\eta)=0$ for all $n$) because we know that $\eta\cdot\nu=0\in \pi_*(S^0)$. The reader can find more property of extensions in \cite{LWX126}.
    \item Rows with \verb|reason=`Syn'|. This means that the differential is obtained by the Generalized Leibniz Rule from the former row. The map used in the Generalized Leibniz Rule and the conditions that the machine checked is written in the \verb|info| column.
    \item Rows with \verb|reason=`SynCs'|. This means that the row is an extension and is obtained by the Generalized Mahowald Trick. The cofiber sequence used in the Generalized Mahowald Trick and the conditions that the machine checked is written in the \verb|info| column.
    \item Rows with \verb|reason=`SynIn'|. This means that by the Generalized Mahowald Trick we determine that an element in the Adams spectral sequence should be a permanent cycle and become a summand of an extension.
    \item Rows with \verb|reason=`TI', `DI', `GI'|. These are essentially the same as rows with \verb|reason=`T', `D', `G'| except that the columns \verb|stem,s,t| refer to the degree of \verb|dx| in stead of \verb|x|. For consecutive rows with \verb|reason=`TI'| the machine are enumerating possible sources of some differential (when we know \verb|dx| is hit by some $d_r$ differential beforehand) instead of the possible targets and the conclusion row is one with \verb|reason=`DI'|.
\end{itemize}

\begin{remark}
    If the reader wants to search for a proof of certain differential, then usually you need to search for the bidegree of the source of the differential with \texttt{reason=`D' or reason=`T'}. Occasionally the machine knows the target should be hit by a certain differential first and deduces the differential by enumerating the possible sources. In this case you also need to search for the bidegree of the target of the differential with \texttt{reason=`DI' or reason=`TI'}.
\end{remark}

\section{Charts and tables}\label{sec:charts}
In this section we list the Adams differentials of $S^0$ in the range $122\le stem\le 127$, $0\le s\le 25$ and Adams differentials of $C\nu$ in the range $stem=126$, $9\le s\le 14$. These differentials would play an important role in our solution of the Last Kervaire Invariant Problem \cite{LWX126}.
\begin{table}
    \centering
\scalebox{1}{\begin{tabular}{|c|l|l|l|}\hline
    $s$ & Elements & $d_r$ & value\\\hline\hline
    \multirow{2}{*}{25} & $e_0g^3\Delta h_1g$ & $d_{2}^{-1}$ & $g^2\Delta^2m$ \\\cline{2-4}
     & $d_0^3g[B_4]$ & $d_{4}$ & $d_0^4x_{65,13}$ \\\hline\hline
    \multirow{1}{*}{24} & $d_0g\Delta^2g^2$ & $d_{2}$ & $d_0e_0g^3\Delta h_2^2$ \\\hline\hline
    \multirow{2}{*}{23} & $Ph_1x_{113,18,2}$ & $d_{3}^{-1}$ & $x_{123,20}$ \\\cline{2-4}
     & $h_0^2d_0x_{108,17}$ & $d_{3}$ & $d_0^4Mg$ \\\hline\hline
    \multirow{3}{*}{22} & $e_0g^2Mg$ & $d_{3}^{-1}$ & $h_0^2h_3x_{116,16}$ \\\cline{2-4}
     & $x_{122,22}$ & $d_{3}$ & $h_0^2d_0x_{107,19}$ \\\cline{2-4}
     & $h_0d_0x_{108,17}$ & $d_{2}$ & $h_0^2d_0x_{107,18}$ \\\hline\hline
    \multirow{2}{*}{21} & $h_0^2d_0e_0x_{91,11}$ & & Permanent \\\cline{2-4}
     & $d_0x_{108,17}$ & $d_{2}$ & $d_0e_0\Delta h_2^2[B_4]+h_0d_0x_{107,18}$ \\\hline\hline
    \multirow{3}{*}{20} & $h_0^4x_{122,16}$ & $d_{2}^{-1}$ & $h_0h_3x_{116,16}$ \\\cline{2-4}
     & $g^3(C_0+h_0^6h_5^2)$ & $d_{3}$ & $g^3\Delta h_2^2n$ \\\cline{2-4}
     & $h_0d_0e_0x_{91,11}$ & $d_{2}$ & $h_0^6x_{121,16}$ \\\hline\hline
    \multirow{2}{*}{19} & $h_0^3x_{122,16}$ & $d_{2}^{-1}$ & $h_3x_{116,16}$ \\\cline{2-4}
     & $d_0e_0x_{91,11}$ & $d_{2}$ & $h_0d_0^2x_{93,12}$ \\\hline\hline
    \multirow{3}{*}{18} & $h_0^2x_{122,16}$ & $d_{3}^{-1}$ & $h_0^2x_{123,13,2}$ \\\cline{2-4}
     & $h_1x_{121,17}$ & $d_{7}^{-1}$ & $x_{123,11,2}+x_{123,11}+h_0h_6[B_4]$ \\\cline{2-4}
     & $d_0x_{108,14}$ & $d_{3}$ & $h_0d_0^2x_{93,12}+h_0^5x_{121,16}$ \\\hline\hline
    \multirow{2}{*}{17} & $h_0x_{122,16}$ & $d_{3}^{-1}$ & $h_0x_{123,13,2}$ \\\cline{2-4}
     & $g^3[H_1]$ & $d_{3}^{-1}$ & $\Delta h_2^2x_{93,8}$ \\\hline\hline
    \multirow{3}{*}{16} & $x_{122,16}+h_0x_{122,15,2}$ & $d_{3}^{-1}$ & $x_{123,13,2}$ \\\cline{2-4}
     & $\Delta h_2^2x_{92,10}$ & $d_{3}^{-1}$ & $x_{123,13}$ \\\cline{2-4}
     & $h_0x_{122,15,2}$ & & Permanent \\\hline\hline
    \multirow{2}{*}{15} & $x_{122,15}$ & $d_{3}$ & $g^3A$ \\\cline{2-4}
     & $x_{122,15,2}$ & $d_{2}$ & $h_0^2h_4x_{106,14}$ \\\hline\hline
    \multirow{1}{*}{14} & $h_0x_{122,13}$ & $d_{3}^{-1}$ & $x_{123,11,2}$ \\\hline\hline
    \multirow{3}{*}{13} & $h_0^2h_6Md_0$ & $d_{2}^{-1}$ & $x_{123,11}$ \\\cline{2-4}
     & $h_1^2x_{120,11}$ & & Permanent \\\cline{2-4}
     & $x_{122,13}$ & $d_{3}$ & $h_0h_4x_{106,14}$ \\\hline\hline
    \multirow{3}{*}{12} & $h_0h_6Md_0$ & $d_{2}^{-1}$ & $x_{123,10}$ \\\cline{2-4}
     & $h_0x_{122,11}$ & $d_{3}^{-1}$ & $h_0x_{123,8}$ \\\cline{2-4}
     & $h_5x_{91,11}$ & & Permanent \\\hline\hline
    \multirow{2}{*}{11} & $x_{122,11}+h_6Md_0$ & $d_{3}^{-1}$ & $x_{123,8}$ \\\cline{2-4}
     & $h_6Md_0$ & & Permanent \\\hline\hline
    \multirow{1}{*}{9-10} & \multicolumn{3}{c|}{}\\\hline\hline
    \multirow{1}{*}{8} & $h_1x_{121,7}$ & $d_{6}$ & $?$ \\\hline\hline
    \multirow{1}{*}{0-7} & \multicolumn{3}{c|}{}\\\hline
\end{tabular}}
\caption{The classical Adams spectral sequence of $S^0$ for $s \le 25$ in stem 122}\vspace{1cm}
\label{Table:S122}
\end{table}
%\newpage

\begin{table}
    \centering
\scalebox{0.95}{\begin{tabular}{|c|l|l|l|}\hline
$s$ & Elements & $d_r$ & value\\\hline\hline
\multirow{1}{*}{25} & \multicolumn{3}{c|}{}\\\hline\hline
\multirow{3}{*}{24} & $h_1Ph_1x_{113,18,2}$ & $d_{2}^{-1}$ & $x_{124,22}$ \\\cline{2-4}
 & $d_0^2\Delta h_2^2Mg$ & $d_{2}^{-1}$ & $d_0x_{110,18}$ \\\cline{2-4}
 & $d_0M\!Px_{56,10}$ & $d_{4}$ & $M\!Px_{69,18}$ \\\hline\hline
\multirow{1}{*}{23} & $g^2\Delta^2m$ & $d_{2}$ & $e_0g^3\Delta h_1g$ \\\hline\hline
\multirow{1}{*}{21-22} & \multicolumn{3}{c|}{}\\\hline\hline
\multirow{1}{*}{20} & $x_{123,20}$ & $d_{3}$ & $Ph_1x_{113,18,2}$ \\\hline\hline
\multirow{4}{*}{19} & $e_0g^2x_{66,7}+h_0^6x_{123,13,2}$ & $d_{2}^{-1}$ & $x_{124,17}$ \\\cline{2-4}
 & $h_0^6x_{123,13,2}$ & $d_{2}^{-1}$ & $h_0^3x_{124,14,2}$ \\\cline{2-4}
 & $h_0e_0x_{106,14}+h_0^2h_3x_{116,16}$ & $d_{3}^{-1}$ & $e_0x_{107,12}$ \\\cline{2-4}
 & $h_0^2h_3x_{116,16}$ & $d_{3}$ & $e_0g^2Mg$ \\\hline\hline
\multirow{3}{*}{18} & $h_0^5x_{123,13,2}$ & $d_{2}^{-1}$ & $h_0^2x_{124,14,2}$ \\\cline{2-4}
 & $e_0x_{106,14}+h_0h_3x_{116,16}$ & & Permanent \\\cline{2-4}
 & $h_0h_3x_{116,16}$ & $d_{2}$ & $h_0^4x_{122,16}$ \\\hline\hline
\multirow{4}{*}{17} & $h_0^2x_{123,15}+h_0^4x_{123,13,2}$ & $d_{2}^{-1}$ & $h_0x_{124,14}$ \\\cline{2-4}
 & $h_0^4x_{123,13,2}$ & $d_{2}^{-1}$ & $h_0x_{124,14,2}+h_0x_{124,14}$ \\\cline{2-4}
 & $d_0x_{109,13}$ & & Permanent \\\cline{2-4}
 & $h_3x_{116,16}$ & $d_{2}$ & $h_0^3x_{122,16}$ \\\hline\hline
\multirow{3}{*}{16} & $h_0x_{123,15}+h_0^3x_{123,13,2}$ & $d_{2}^{-1}$ & $x_{124,14}$ \\\cline{2-4}
 & $h_0^3x_{123,13,2}$ & $d_{2}^{-1}$ & $x_{124,14,2}+x_{124,14}$ \\\cline{2-4}
 & $h_1x_{122,15,2}$ & $d_{3}^{-1}$ & $h_4x_{109,12}$ \\\hline\hline
\multirow{3}{*}{15} & $x_{123,15}$ & $d_{4}^{-1}$ & $x_{124,11,2}+x_{124,11}$ \\\cline{2-4}
 & $h_4x_{108,14}$ & $d_{5}^{-1}$ & $h_0x_{124,9}$ \\\cline{2-4}
 & $h_0^2x_{123,13,2}$ & $d_{3}$ & $h_0^2x_{122,16}$ \\\hline\hline
\multirow{3}{*}{14} & $h_0^3x_{123,11}$ & $d_{2}^{-1}$ & $h_0x_{124,11}$ \\\cline{2-4}
 & $h_0x_{123,13,2}$ & $d_{3}$ & $h_0x_{122,16}$ \\\cline{2-4}
 & $\Delta h_2^2x_{93,8}$ & $d_{3}$ & $g^3[H_1]$ \\\hline\hline
\multirow{3}{*}{13} & $h_0^2x_{123,11}$ & $d_{2}^{-1}$ & $x_{124,11}$ \\\cline{2-4}
 & $x_{123,13,2}$ & $d_{3}$ & $x_{122,16}+h_0x_{122,15,2}$ \\\cline{2-4}
 & $x_{123,13}$ & $d_{3}$ & $\Delta h_2^2x_{92,10}$ \\\hline\hline
\multirow{3}{*}{12} & $h_0x_{123,11}+h_0^2h_6[B_4]$ & $d_{3}^{-1}$ & $x_{124,9,2}+h_0x_{124,8}$ \\\cline{2-4}
 & $x_{123,12}$ & $d_{3}^{-1}$ & $x_{124,9}+h_0x_{124,8}$ \\\cline{2-4}
 & $h_0^2h_6[B_4]$ & $d_{5}^{-1}$ & $h_6A$ \\\hline\hline
\multirow{5}{*}{11} & $h_0^2x_{123,9}$ & $d_{2}^{-1}$ & $x_{124,9}$ \\\cline{2-4}
 & $h_5x_{92,10}$ & $d_{5}^{-1}$ & $x_{124,6}$ \\\cline{2-4}
 & $x_{123,11,2}+x_{123,11}+h_0h_6[B_4]$ & $d_{7}$ & $h_1x_{121,17}$ \\\cline{2-4}
 & $x_{123,11}+h_0h_6[B_4]$ & $d_{3}$ & $h_0x_{122,13}$ \\\cline{2-4}
 & $h_0h_6[B_4]$ & $d_{2}$ & $h_0^2h_6Md_0$ \\\hline\hline
\multirow{3}{*}{10} & $h_0x_{123,9}$ & $d_{2}^{-1}$ & $x_{124,8}$ \\\cline{2-4}
 & $x_{123,10}+h_6[B_4]$ & $d_{3}^{-1}$ & $x_{124,7}$ \\\cline{2-4}
 & $h_6[B_4]$ & $d_{2}$ & $h_0h_6Md_0$ \\\hline\hline
\multirow{2}{*}{9} & $x_{123,9}+h_0x_{123,8}$ & $d_{12}$ & $?$ \\\cline{2-4}
 & $h_0x_{123,8}$ & $d_{3}$ & $h_0x_{122,11}+h_0h_6Md_0$ \\\hline\hline
\multirow{1}{*}{8} & $x_{123,8}$ & $d_{3}$ & $x_{122,11}+h_6Md_0$ \\\hline\hline
\multirow{1}{*}{0-7} & \multicolumn{3}{c|}{}\\\hline
\end{tabular}}
\caption{The classical Adams spectral sequence of $S^0$ for $s \le 25$ in stem 123}
\label{Table:S123}
\end{table}

%\newpage

\begin{table}
    \centering
\scalebox{0.95}{\begin{tabular}{|c|l|l|l|}\hline
$s$ & Elements & $d_r$ & value\\\hline\hline
\multirow{3}{*}{25} & $h_0^{11}x_{124,14,2}$ & $d_{2}^{-1}$ & $h_0^9x_{125,14}$ \\\cline{2-4}
 & $ix_{101,18}$ & $d_{2}$ & $d_0^2\Delta h_2^2x_{65,13}+h_0Pd_0x_{101,18}$ \\\cline{2-4}
 & $d_0e_0\Delta^3h_1g$ & $d_{2}$ & $d_0^2g^3m$ \\\hline\hline
\multirow{2}{*}{24} & $h_0^{10}x_{124,14,2}$ & $d_{2}^{-1}$ & $h_0^8x_{125,14}$ \\\cline{2-4}
 & $h_0^2d_0x_{110,18}$ & $d_{2}^{-1}$ & $h_0e_0x_{108,17}$ \\\hline\hline
\multirow{4}{*}{23} & $h_0^9x_{124,14,2}$ & $d_{2}^{-1}$ & $h_0^7x_{125,14}$ \\\cline{2-4}
 & $d_0g\Delta h_2^2[B_4]+h_0d_0x_{110,18}$ & $d_{2}^{-1}$ & $e_0x_{108,17}$ \\\cline{2-4}
 & $h_0d_0x_{110,18}$ & $d_{3}^{-1}$ & $x_{125,20}$ \\\cline{2-4}
 & $d_0x_{110,19}$ & & Permanent \\\hline\hline
\multirow{4}{*}{22} & $h_0^8x_{124,14,2}$ & $d_{2}^{-1}$ & $h_0^6x_{125,14}$ \\\cline{2-4}
 & $g^2\Delta^2t$ & $d_{3}^{-1}$ & $gx_{105,15}$ \\\cline{2-4}
 & $x_{124,22}$ & $d_{2}$ & $h_1Ph_1x_{113,18,2}$ \\\cline{2-4}
 & $d_0x_{110,18}$ & $d_{2}$ & $d_0^2\Delta h_2^2Mg$ \\\hline\hline
\multirow{2}{*}{21} & $h_0^7x_{124,14,2}$ & $d_{2}^{-1}$ & $h_0^5x_{125,14}$ \\\cline{2-4}
 & $h_0d_0^2[\Delta\Delta_1g]$ & $d_{2}^{-1}$ & $d_0gx_{91,11}$ \\\hline\hline
\multirow{2}{*}{20} & $h_0^6x_{124,14,2}$ & $d_{2}^{-1}$ & $h_0^4x_{125,14}$ \\\cline{2-4}
 & $d_0^2[\Delta\Delta_1g]$ & $d_{5}^{-1}$ & $h_1x_{124,14}$ \\\hline\hline
\multirow{2}{*}{19} & $h_0^5x_{124,14,2}$ & $d_{2}^{-1}$ & $h_0^3x_{125,14}$ \\\cline{2-4}
 & $d_0x_{110,15}$ & & Permanent \\\hline\hline
\multirow{2}{*}{18} & $h_0x_{124,17}+h_0^4x_{124,14,2}$ & $d_{2}^{-1}$ & $x_{125,16}$ \\\cline{2-4}
 & $h_0^4x_{124,14,2}$ & $d_{2}^{-1}$ & $h_0^2x_{125,14}$ \\\hline\hline
\multirow{4}{*}{17} & $h_0^3x_{124,14}$ & $d_{2}^{-1}$ & $x_{125,15}$ \\\cline{2-4}
 & $h_0^2x_{124,15}$ & $d_{2}^{-1}$ & $h_0x_{125,14}$ \\\cline{2-4}
 & $x_{124,17}$ & $d_{2}$ & $e_0g^2x_{66,7}+h_0^6x_{123,13,2}$ \\\cline{2-4}
 & $h_0^3x_{124,14,2}$ & $d_{2}$ & $h_0^6x_{123,13,2}$ \\\hline\hline
\multirow{5}{*}{16} & $h_0x_{124,15}$ & $d_{2}^{-1}$ & $x_{125,14}$ \\\cline{2-4}
 & $h_1x_{123,15}$ & $d_{3}^{-1}$ & $h_3x_{118,12}$ \\\cline{2-4}
 & $h_0^2x_{124,14}$ & & Permanent \\\cline{2-4}
 & $e_0x_{107,12}$ & $d_{3}$ & $e_0g^2x_{66,7}+h_0e_0x_{106,14}+h_0^2h_3x_{116,16}$ \\\cline{2-4}
 & $h_0^2x_{124,14,2}$ & $d_{2}$ & $h_0^5x_{123,13,2}$ \\\hline\hline
\multirow{4}{*}{15} & $x_{124,15}$ & $d_{4}^{-1}$ & $h_6x_{62,10}$ \\\cline{2-4}
 & $h_3^2x_{110,13}+h_0x_{124,14}$ & & Permanent \\\cline{2-4}
 & $h_0x_{124,14}$ & $d_{2}$ & $h_0^2x_{123,15}+h_0^4x_{123,13,2}$ \\\cline{2-4}
 & $h_0x_{124,14,2}$ & $d_{2}$ & $h_0^2x_{123,15}$ \\\hline\hline
\multirow{5}{*}{14} & $h_1x_{123,13}$ & $d_{2}^{-1}$ & $x_{125,12}$ \\\cline{2-4}
 & $h_1x_{123,13,2}$ & $d_{2}^{-1}$ & $x_{125,12,2}$ \\\cline{2-4}
 & $\Delta h_2^2x_{94,8}$ & & Permanent \\\cline{2-4}
 & $x_{124,14}$ & $d_{2}$ & $h_0x_{123,15}+h_0^3x_{123,13,2}$ \\\cline{2-4}
 & $x_{124,14,2}$ & $d_{2}$ & $h_0x_{123,15}$ \\\hline\hline
\multirow{4}{*}{13} & $h_0^5x_{124,8}$ & $d_{2}^{-1}$ & $h_0^3x_{125,8}$ \\\cline{2-4}
 & $[H_1](\Delta e_1+C_0+h_0^6h_5^2)$ & $d_{3}^{-1}$ & $x_{125,10,2}$ \\\cline{2-4}
 & $e_0\Delta h_6g$ & & Permanent \\\cline{2-4}
 & $h_4x_{109,12}$ & $d_{3}$ & $h_1x_{122,15,2}$ \\\hline
\end{tabular}}
\caption{The classical Adams spectral sequence of $S^0$ for $13 \le s \le 25$ in stem 124}
\label{Table:S124.13}
\end{table}

%\newpage

\begin{table}
    \centering
\scalebox{0.9}{\begin{tabular}{|c|l|l|l|}\hline %%%%%%%%%%%%%%%%%%%%%%
\multirow{5}{*}{12} & $h_0x_{124,11,2}+h_0x_{124,11}$ & $d_{2}^{-1}$ & $x_{125,10}$ \\\cline{2-4}
 & $h_0^2x_{124,10,2}+h_0^4x_{124,8}$ & $d_{2}^{-1}$ & $h_0x_{125,9,2}$ \\\cline{2-4}
 & $h_0^4x_{124,8}$ & $d_{2}^{-1}$ & $h_0^2x_{125,8}$ \\\cline{2-4}
 & $h_1x_{123,11,2}$ & $d_{3}^{-1}$ & $x_{125,9}$ \\\cline{2-4}
 & $h_0x_{124,11}$ & $d_{2}$ & $h_0^3x_{123,11}$ \\\hline\hline
\multirow{5}{*}{11} & $h_0x_{124,10,2}+h_0^3x_{124,8}$ & $d_{2}^{-1}$ & $x_{125,9,2}$ \\\cline{2-4}
 & $h_0^3x_{124,8}$ & $d_{2}^{-1}$ & $h_0x_{125,8}$ \\\cline{2-4}
 & $x_{124,11,3}$ & $d_{3}^{-1}$ & $x_{125,8,2}$ \\\cline{2-4}
 & $x_{124,11,2}+x_{124,11}$ & $d_{4}$ & $x_{123,15}$ \\\cline{2-4}
 & $x_{124,11}$ & $d_{2}$ & $h_0^2x_{123,11}$ \\\hline\hline
\multirow{5}{*}{10} & $h_1x_{123,9}+h_0^2x_{124,8}$ & $d_{2}^{-1}$ & $x_{125,8}$ \\\cline{2-4}
 & $x_{124,10,2}+h_0x_{124,9}$ & $d_{4}^{-1}$ & $h_6[H_1]$ \\\cline{2-4}
 & $x_{124,10}+h_0^2x_{124,8}$ & $d_{4}^{-1}$ & $h_6[H_1]+h_0x_{125,5}$ \\\cline{2-4}
 & $h_0^2x_{124,8}$ & & Permanent \\\cline{2-4}
 & $h_0x_{124,9}$ & $d_{5}$ & $h_4x_{108,14}$ \\\hline\hline
\multirow{3}{*}{9} & $x_{124,9,2}+h_0x_{124,8}$ & $d_{3}$ & $h_0x_{123,11}+h_0^2h_6[B_4]$ \\\cline{2-4}
 & $x_{124,9}+h_0x_{124,8}$ & $d_{3}$ & $x_{123,12}$ \\\cline{2-4}
 & $h_0x_{124,8}$ & $d_{2}$ & $h_0^2x_{123,9}$ \\\hline\hline
\multirow{1}{*}{8} & $x_{124,8}$ & $d_{2}$ & $h_0x_{123,9}$ \\\hline\hline
\multirow{3}{*}{7} & $h_0x_{124,6}$ & $d_{2}^{-1}$ & $x_{125,5}$ \\\cline{2-4}
 & $h_6A$ & $d_{5}$ & $h_0^2h_6[B_4]$ \\\cline{2-4}
 & $x_{124,7}$ & $d_{3}$ & $x_{123,10}+h_6[B_4]$ \\\hline\hline
\multirow{1}{*}{6} & $x_{124,6}$ & $d_{5}$ & $h_5x_{92,10}$ \\\hline\hline
\multirow{1}{*}{0-5} & \multicolumn{3}{c|}{}\\\hline
\end{tabular}}
\caption{The classical Adams spectral sequence of $S^0$ for $s \le 12$ in stem 124}
\label{Table:S124.12}
\end{table}

\begin{table}
    \centering
\scalebox{0.9}{\begin{tabular}{|c|l|l|l|}\hline
$s$ & Elements & $d_r$ & value\\\hline\hline
\multirow{4}{*}{25} & $d_0^2e_0g[B_4]$ & $d_{4}^{-1}$ & $h_1x_{125,20}$ \\\cline{2-4}
 & $x_{125,25,2}+x_{125,25}+g^4\Delta h_1g$ & $d_{4}^{-1}$ & $x_{126,21}+\text{possibly }h_1x_{125,20}$ \\\cline{2-4}
 & $g^4\Delta h_1g$ & & Permanent \\\cline{2-4}
 & $x_{125,25}$ & $d_{2}$ & $h_0x_{124,26}$ \\\hline\hline
\multirow{1}{*}{24} & $e_0g\Delta^2g^2$ & $d_{2}$ & $d_0g^4\Delta h_2^2$ \\\hline\hline
\multirow{2}{*}{23} & $h_0^2e_0x_{108,17}$ & $d_{3}$ & $d_0^3e_0Mg$ \\\cline{2-4}
 & $h_0^9x_{125,14}$ & $d_{2}$ & $h_0^{11}x_{124,14,2}$ \\\hline\hline
\multirow{4}{*}{22} & $g^3Mg$ & $d_{4}^{-1}$ & $x_{126,18}$ \\\cline{2-4}
 & $ix_{102,15}+h_0^8x_{125,14}$ & $d_{4}^{-1}$ & $gx_{106,14}+e_0x_{109,14,2}$ \\\cline{2-4}
 & $h_0^8x_{125,14}$ & $d_{2}$ & $h_0^{10}x_{124,14,2}$ \\\cline{2-4}
 & $h_0e_0x_{108,17}$ & $d_{2}$ & $h_0^2d_0x_{110,18}$ \\\hline\hline
\multirow{3}{*}{21} & $x_{125,21}$ & $d_{4}^{-1}$ & $d_0x_{112,13}$ \\\cline{2-4}
 & $h_0^7x_{125,14}$ & $d_{2}$ & $h_0^9x_{124,14,2}$ \\\cline{2-4}
 & $e_0x_{108,17}$ & $d_{2}$ & $d_0g\Delta h_2^2[B_4]+h_0d_0x_{110,18}$ \\\hline\hline
\multirow{3}{*}{20} & $h_0d_0gx_{91,11}$ & $d_{2}^{-1}$ & $x_{126,18,2}$ \\\cline{2-4}
 & $x_{125,20}$ & $d_{3}$ & $d_0g\Delta h_2^2[B_4]$ \\\cline{2-4}
 & $h_0^6x_{125,14}$ & $d_{2}$ & $h_0^8x_{124,14,2}$ \\\hline
 \end{tabular}}
\caption{The classical Adams spectral sequence of $S^0$ for $20 \le s \le 25$ in stem 125}
\label{Table:S125.20}
\end{table}

\begin{table}
    \centering
\scalebox{0.85}{\begin{tabular}{|c|l|l|l|}\hline %%%%%%%%%%%%%%%%%%%%%%
\multirow{3}{*}{19} & $gx_{105,15}$ & $d_{3}$ & $g^2\Delta^2t$ \\\cline{2-4}
 & $h_0^5x_{125,14}$ & $d_{2}$ & $h_0^7x_{124,14,2}$ \\\cline{2-4}
 & $d_0gx_{91,11}$ & $d_{2}$ & $h_0d_0^2[\Delta\Delta_1g]$ \\\hline\hline
\multirow{2}{*}{18} & $d_0^2x_{97,10}$ & $d_{5}^{-1}$ & $h_1x_{125,12,2}$ \\\cline{2-4}
 & $h_0^4x_{125,14}$ & $d_{2}$ & $h_0^6x_{124,14,2}$ \\\hline\hline
\multirow{2}{*}{17} & $h_0^2Q_2x_{68,8}$ & $d_{2}^{-1}$ & $h_0D_2x_{68,8}$ \\\cline{2-4}
 & $h_0^3x_{125,14}$ & $d_{2}$ & $h_0^5x_{124,14,2}$ \\\hline\hline
\multirow{4}{*}{16} & $h_0Q_2x_{68,8}$ & $d_{2}^{-1}$ & $D_2x_{68,8}$ \\\cline{2-4}
 & $h_1x_{124,15}$ & $d_{4}^{-1}$ & $h_0^2x_{126,10}$ \\\cline{2-4}
 & $x_{125,16}$ & $d_{2}$ & $h_0x_{124,17}+h_0^4x_{124,14,2}$ \\\cline{2-4}
 & $h_0^2x_{125,14}$ & $d_{2}$ & $h_0^4x_{124,14,2}$ \\\hline\hline
\multirow{5}{*}{15} & $h_0^3x_{125,12}$ & $d_{2}^{-1}$ & $h_0h_3x_{119,11}$ \\\cline{2-4}
 & $Q_2x_{68,8}$ & $d_{4}^{-1}$ & $h_1x_{125,10}$ \\\cline{2-4}
 & $h_1x_{124,14}$ & $d_{5}$ & $d_0^2[\Delta\Delta_1g]$ \\\cline{2-4}
 & $x_{125,15}$ & $d_{2}$ & $h_0^3x_{124,14}$ \\\cline{2-4}
 & $h_0x_{125,14}$ & $d_{2}$ & $h_0^2x_{124,15}$ \\\hline\hline
\multirow{3}{*}{14} & $h_0^2x_{125,12}$ & $d_{2}^{-1}$ & $h_3x_{119,11}$ \\\cline{2-4}
 & $h_1h_4x_{109,12}$ & & Permanent \\\cline{2-4}
 & $x_{125,14}$ & $d_{2}$ & $h_0x_{124,15}$ \\\hline\hline
\multirow{5}{*}{13} & $h_0^4x_{125,9,2}$ & $d_{2}^{-1}$ & $x_{126,11}$ \\\cline{2-4}
 & $h_0^5x_{125,8}$ & $d_{2}^{-1}$ & $h_0x_{126,10}$ \\\cline{2-4}
 & $nx_{94,8}$ & $d_{4}^{-1}$ & $h_1x_{125,8}$ \\\cline{2-4}
 & $h_0x_{125,12}$ & $d_{4}^{-1}$ & $h_1x_{125,8,2}$ \\\cline{2-4}
 & $h_3x_{118,12}$ & $d_{3}$ & $h_1x_{123,15}$ \\\hline\hline
\multirow{5}{*}{12} & $h_0h_6x_{62,10}$ & $d_{2}^{-1}$ & $x_{126,10}$ \\\cline{2-4}
 & $h_0^3x_{125,9,2}+h_0^4x_{125,8}$ & $d_{3}^{-1}$ & $x_{126,9}$ \\\cline{2-4}
 & $h_0^4x_{125,8}$ & $d_{3}^{-1}$ & $x_{126,9}+h_0x_{126,8,3}$ \\\cline{2-4}
 & $x_{125,12}$ & $d_{2}$ & $h_1x_{123,13}$ \\\cline{2-4}
 & $x_{125,12,2}$ & $d_{2}$ & $h_1x_{123,13,2}$ \\\hline\hline
\multirow{5}{*}{11} & $h_1x_{124,10,2}$ & $d_{3}^{-1}$ & $x_{126,8}$ \\\cline{2-4}
 & $h_1x_{124,10}$ & $d_{3}^{-1}$ & $x_{126,8,2}$ \\\cline{2-4}
 & $h_0^2x_{125,9,2}$ & $d_{5}$ & $?$ \\\cline{2-4}
 & $h_6x_{62,10}$ & $d_{4}$ & $x_{124,15}$ \\\cline{2-4}
 & $h_0^3x_{125,8}$ & $d_{2}$ & $h_0^5x_{124,8}$ \\\hline\hline
\multirow{5}{*}{10} & $h_0x_{125,9}$ & $d_{2}^{-1}$ & $x_{126,8,3}$ \\\cline{2-4}
 & $x_{125,10,2}$ & $d_{3}$ & $[H_1](\Delta e_1+C_0+h_0^6h_5^2)$ \\\cline{2-4}
 & $x_{125,10}$ & $d_{2}$ & $h_0x_{124,11,2}+h_0x_{124,11}$ \\\cline{2-4}
 & $h_0x_{125,9,2}$ & $d_{2}$ & $h_0^2x_{124,10,2}+h_0^4x_{124,8}$ \\\cline{2-4}
 & $h_0^2x_{125,8}$ & $d_{2}$ & $h_0^4x_{124,8}$ \\\hline\hline
\multirow{5}{*}{9} & $h_6(\Delta e_1+C_0+h_0^6h_5^2)$ & & Permanent \\\cline{2-4}
 & $h_5x_{94,8}$ & $d_{7}$ & $?$ \\\cline{2-4}
 & $x_{125,9}$ & $d_{3}$ & $h_1x_{123,11,2}$ \\\cline{2-4}
 & $x_{125,9,2}$ & $d_{2}$ & $h_0x_{124,10,2}+h_0^3x_{124,8}$ \\\cline{2-4}
 & $h_0x_{125,8}$ & $d_{2}$ & $h_0^3x_{124,8}$ \\\hline\hline
\multirow{2}{*}{8} & $x_{125,8,2}$ & $d_{3}$ & $x_{124,11,3}$ \\\cline{2-4}
 & $x_{125,8}$ & $d_{2}$ & $h_1x_{123,9}+h_0^2x_{124,8}$ \\\hline\hline
\multirow{1}{*}{7} & $h_0^2x_{125,5}$ & $d_{3}^{-1}$ & $x_{126,4}$ \\\hline\hline
\multirow{2}{*}{6} & $h_6[H_1]$ & $d_{4}$ & $x_{124,10,2}+h_0x_{124,9}$ \\\cline{2-4}
 & $h_0x_{125,5}$ & $d_{4}$ & $x_{124,10,2}+x_{124,10}+h_0x_{124,9}+h_0^2x_{124,8}$ \\\hline\hline
\multirow{1}{*}{5} & $x_{125,5}$ & $d_{2}$ & $h_0x_{124,6}$ \\\hline
\hline
\multirow{1}{*}{0-4} & \multicolumn{3}{c|}{}\\\hline
 \end{tabular}}
\caption{The classical Adams spectral sequence of $S^0$ for $s \le 19$ in stem 125}
\label{Table:S125.19}
\end{table}

\begin{table}
    \centering
\scalebox{0.85}{\begin{tabular}{|c|l|l|l|}\hline
$s$ & Elements & $d_r$ & value\\\hline\hline
\multirow{1}{*}{25} & $h_0^7x_{126,18}$ & $d_{3}^{-1}$ & $h_0^{21}h_7$ \\\hline\hline
\multirow{4}{*}{24} & $d_0e_0\Delta h_2^2Mg$ & $d_{2}^{-1}$ & $d_0x_{113,18}$ \\\cline{2-4}
 & $h_0^6x_{126,18}$ & $d_{3}^{-1}$ & $h_0^{20}h_7$ \\\cline{2-4}
 & $g^4\Delta h_2c_1$ & $d_{3}^{-1}$ & $g^3C^{\prime\prime}$ \\\cline{2-4}
 & $d_0Pd_0M^2$ & $d_{4}^{-1}$ & $d_0e_0[\Delta\Delta_1g]$ \\\hline\hline
\multirow{2}{*}{23} & $h_0^5x_{126,18}$ & $d_{3}^{-1}$ & $h_0^{19}h_7$ \\\cline{2-4}
 & $x_{126,23}$ & $d_{4}^{-1}$ & $e_0x_{110,15}$ \\\hline\hline
\multirow{2}{*}{22} & $h_0x_{126,21}+h_0^4x_{126,18}$ & $d_{3}^{-1}$ & $h_1x_{126,18,2}$ \\\cline{2-4}
 & $h_0^4x_{126,18}$ & $d_{3}^{-1}$ & $h_0^{18}h_7$ \\\hline\hline
\multirow{3}{*}{21} & $h_0^3x_{126,18}$ & $d_{3}^{-1}$ & $h_0^{17}h_7$ \\\cline{2-4}
 & $h_1x_{125,20}$ & $d_{4}$ & $d_0^2e_0g[B_4]$ \\\cline{2-4}
 & $x_{126,21}$ & $d_{4}$ & $x_{125,25,2}+x_{125,25}+g^4\Delta h_1g+\text{possibly }d_0^2e_0gB_4$ \\\hline\hline
\multirow{2}{*}{20} & $h_0^2x_{126,18}$ & $d_{3}^{-1}$ & $h_0^{16}h_7$ \\\cline{2-4}
 & $d_0x_{112,16}$ & $d_{5}^{-1}$ & $x_{127,15}$ \\\hline\hline
\multirow{2}{*}{19} & $h_0x_{126,18}$ & $d_{3}^{-1}$ & $h_0^{15}h_7$ \\\cline{2-4}
 & $g^3x_{66,7}$ & $d_{3}^{-1}$ & $gx_{107,12}$ \\\hline\hline
\multirow{4}{*}{18} & $x_{126,18}+e_0x_{109,14,2}$ & $d_{7}$ & $?$ \\\cline{2-4}
 & $e_0x_{109,14,2}$ & $d_{4}$ & $g^3Mg$ \\\cline{2-4}
 & $gx_{106,14}$ & $d_{4}$ & $ix_{102,15}+g^3Mg+h_0^8x_{125,14}$ \\\cline{2-4}
 & $x_{126,18,2}$ & $d_{2}$ & $h_0d_0gx_{91,11}$ \\\hline\hline
\multirow{4}{*}{17} & $h_0^{15}h_6^2$ & $d_{2}^{-1}$ & $h_0^{14}h_7$ \\\cline{2-4}
 & $h_1^2x_{124,15}$ & $d_{2}^{-1}$ & $h_6x_{64,14}$ \\\cline{2-4}
 & $x_{126,17}$ & $d_{8}$ & $?$ \\\cline{2-4}
 & $d_0x_{112,13}$ & $d_{4}$ & $x_{125,21}$ \\\hline\hline
\multirow{3}{*}{16} & $h_0^{14}h_6^2$ & $d_{2}^{-1}$ & $h_0^{13}h_7$ \\\cline{2-4}
 & $h_0^2D_2x_{68,8}$ & $d_{3}^{-1}$ & $x_{127,13}$ \\\cline{2-4}
 & $h_1^2x_{124,14}$ & $d_{6}^{-1}$ & $h_2x_{124,9}+h_0^2x_{127,8}$ \\\hline\hline
\multirow{2}{*}{15} & $h_0^{13}h_6^2$ & $d_{2}^{-1}$ & $h_0^{12}h_7$ \\\cline{2-4}
 & $h_0D_2x_{68,8}$ & $d_{2}$ & $h_0^2Q_2x_{68,8}$ \\\hline\hline
\multirow{4}{*}{14} & $h_0^{12}h_6^2$ & $d_{2}^{-1}$ & $h_0^{11}h_7$ \\\cline{2-4}
 & $x_{126,14}$ & $d_{4}^{-1}$ & $h_0^2x_{127,8}$ \\\cline{2-4}
 & $h_1h_3x_{118,12}$ & $d_{5}^{-1}$ & $h_1x_{126,8,2}$ \\\cline{2-4}
 & $D_2x_{68,8}$ & $d_{2}$ & $h_0Q_2x_{68,8}$ \\\hline\hline
\multirow{3}{*}{13} & $h_0^{11}h_6^2$ & $d_{2}^{-1}$ & $h_0^{10}h_7$ \\\cline{2-4}
 & $h_1x_{125,12,2}$ & $d_{5}$ & $d_0^2x_{97,10}$ \\\cline{2-4}
 & $h_0h_3x_{119,11}$ & $d_{2}$ & $h_0^3x_{125,12}$ \\\hline\hline
\multirow{5}{*}{12} & $d_1x_{94,8}$ & $d_{2}^{-1}$ & $x_{127,10}$ \\\cline{2-4}
 & $h_0x_{126,11}$ & $d_{2}^{-1}$ & $h_3x_{120,9}$ \\\cline{2-4}
 & $h_0^{10}h_6^2$ & $d_{2}^{-1}$ & $h_0^9h_7$ \\\cline{2-4}
 & $h_0^2x_{126,10}$ & $d_{4}$ & $h_1x_{124,15}$ \\\cline{2-4}
 & $h_3x_{119,11}$ & $d_{2}$ & $h_0^2x_{125,12}$ \\\hline\hline
\multirow{6}{*}{11} & $h_0^2x_{126,9}$ & $d_{2}^{-1}$ & $h_0x_{127,8}$ \\\cline{2-4}
 & $h_0^9h_6^2$ & $d_{2}^{-1}$ & $h_0^8h_7$ \\\cline{2-4}
 & $h_1x_{125,10,2}+h_1x_{125,10}$ & & Permanent \\\cline{2-4}
 & $h_1x_{125,10}$ & $d_{4}$ & $Q_2x_{68,8}$ \\\cline{2-4}
 & $x_{126,11}$ & $d_{2}$ & $h_0^4x_{125,9,2}$ \\\cline{2-4}
 & $h_0x_{126,10}$ & $d_{2}$ & $h_0^5x_{125,8}$ \\\hline
 \end{tabular}}
\caption{The classical Adams spectral sequence of $S^0$ for $11 \le s \le 25$ in stem 126}
\label{Table:S126.11}
\end{table}

\begin{table}
    \centering
\scalebox{0.85}{\begin{tabular}{|c|l|l|l|}\hline %%%%%%%%%%%%%%%%%%%%%%
\multirow{5}{*}{10} & $h_0x_{126,9}$ & $d_{2}^{-1}$ & $x_{127,8}$ \\\cline{2-4}
 & $h_0^2x_{126,8}$ & $d_{2}^{-1}$ & $h_0x_{127,7}$ \\\cline{2-4}
 & $h_0^8h_6^2$ & $d_{2}^{-1}$ & $h_0^7h_7$ \\\cline{2-4}
 & $h_0^2x_{126,8,3}$ & & Permanent \\\cline{2-4}
 & $x_{126,10}$ & $d_{2}$ & $h_0h_6x_{62,10}$ \\\hline\hline
\multirow{6}{*}{9} & $h_0x_{126,8}$ & $d_{2}^{-1}$ & $x_{127,7}$ \\\cline{2-4}
 & $h_0^7h_6^2$ & $d_{2}^{-1}$ & $h_0^6h_7$ \\\cline{2-4}
 & $h_1x_{125,8}$ & $d_{4}$ & $nx_{94,8}$ \\\cline{2-4}
 & $h_1x_{125,8,2}$ & $d_{4}$ & $h_0x_{125,12}$ \\\cline{2-4}
 & $x_{126,9}$ & $d_{3}$ & $h_0^3x_{125,9,2}+h_0^4x_{125,8}$ \\\cline{2-4}
 & $h_0x_{126,8,3}$ & $d_{3}$ & $h_0^3x_{125,9,2}$ \\\hline\hline
\multirow{6}{*}{8} & $h_0^6h_6^2$ & $d_{2}^{-1}$ & $h_0^5h_7$ \\\cline{2-4}
 & $h_6(C^{\prime}+X_2)$ & $d_{17}$ & $?$ \\\cline{2-4}
 & $x_{126,8,4}+x_{126,8}$ & $d_{6}$ & $?$ \\\cline{2-4}
 & $x_{126,8}$ & $d_{3}$ & $h_1x_{124,10,2}$ \\\cline{2-4}
 & $x_{126,8,2}$ & $d_{3}$ & $h_1x_{124,10}$ \\\cline{2-4}
 & $x_{126,8,3}$ & $d_{2}$ & $h_0x_{125,9}$ \\\hline\hline
\multirow{2}{*}{7} & $h_0^5h_6^2$ & $d_{2}^{-1}$ & $h_0^4h_7$ \\\cline{2-4}
 & $h_1h_6[H_1]$ & $d_{18}$ & $?$ \\\hline\hline
\multirow{2}{*}{6} & $h_0^4h_6^2$ & $d_{2}^{-1}$ & $h_0^3h_7$ \\\cline{2-4}
 & $x_{126,6}$ & $d_{3}$ & $h_5x_{94,8} +\text{possibly }h_6(\Delta e_1 + C_0 + h_0^6h_5^2)$ \\\hline\hline
\multirow{1}{*}{5} & $h_0^3h_6^2$ & $d_{2}^{-1}$ & $h_0^2h_7$ \\\hline\hline
\multirow{2}{*}{4} & $h_0^2h_6^2$ & $d_{2}^{-1}$ & $h_0h_7$ \\\cline{2-4}
 & $x_{126,4}$ & $d_{3}$ & $h_0^2x_{125,5}$ \\\hline\hline
\multirow{1}{*}{3} & $h_0h_6^2$ & $d_{2}^{-1}$ & $h_7$ \\\hline\hline
\multirow{1}{*}{2} & $h_6^2$ & $d_{7}$ & $?$ \\\hline\hline
\multirow{1}{*}{0-1} & \multicolumn{3}{c|}{}\\\hline
\end{tabular}}
\caption{The classical Adams spectral sequence of $S^0$ for $s \le 10$ in stem 126}
\label{Table:S126.10}
\end{table}

\begin{table}
    \centering
\scalebox{0.85}{\begin{tabular}{|c|l|l|l|}\hline
$s$ & Elements & $d_r$ & value\\\hline\hline
\multirow{3}{*}{25} & $h_0^{24}h_7$ & $d_{3}$ & $h_0^{10}x_{126,18}$ \\\cline{2-4}
 & $ix_{104,18}$ & $d_{2}$ & $d_0^3x_{84,15,2}+h_0d_0x_{112,22}$ \\\cline{2-4}
 & $d_0g\Delta^3h_1g$ & $d_{2}$ & $d_0e_0g^3m$ \\\hline\hline
\multirow{3}{*}{24} & $h_0^2d_0x_{113,18}$ & $d_{2}^{-1}$ & $h_0gx_{108,17}$ \\\cline{2-4}
 & $h_1x_{126,23}$ & $d_{3}^{-1}$ & $x_{128,21}$ \\\cline{2-4}
 & $h_0^{23}h_7$ & $d_{3}$ & $h_0^9x_{126,18}$ \\\hline\hline
\multirow{4}{*}{23} & $e_0g\Delta h_2^2[B_4]+h_0d_0x_{113,18}$ & $d_{2}^{-1}$ & $gx_{108,17}$ \\\cline{2-4}
 & $h_0d_0x_{113,18}$ & $d_{4}$ & $d_0^3x_{84,15,2}$ \\\cline{2-4}
 & $h_0^{22}h_7$ & $d_{3}$ & $h_0^8x_{126,18}$ \\\cline{2-4}
 & $d_0Pd_0x_{91,11}$ & $d_{3}$ & $h_1x_{125,25}$ \\\hline\hline
\multirow{3}{*}{22} & $d_0x_{113,18,2}$ & $d_{4}^{-1}$ & $d_0e_0x_{97,10}$ \\\cline{2-4}
 & $h_0^{21}h_7$ & $d_{3}$ & $h_0^7x_{126,18}$ \\\cline{2-4}
 & $d_0x_{113,18}$ & $d_{2}$ & $d_0e_0\Delta h_2^2Mg$ \\\hline\hline
\multirow{4}{*}{21} & $h_0^6x_{127,15}$ & $d_{2}^{-1}$ & $h_0^5x_{128,14}$ \\\cline{2-4}
 & $x_{127,21}+g^3C^{\prime\prime}$ & & Permanent \\\cline{2-4}
 & $h_0^{20}h_7$ & $d_{3}$ & $h_0^6x_{126,18}$ \\\cline{2-4}
 & $g^3C^{\prime\prime}$ & $d_{3}$ & $g^4\Delta h_2c_1$ \\\hline
 \end{tabular}}
\caption{The classical Adams spectral sequence of $S^0$ for $21 \le s \le 25$ in stem 127}
\label{Table:S127.21}
\end{table}

\begin{table}
    \centering
\scalebox{0.9}{\begin{tabular}{|c|l|l|l|}\hline %%%%%%%%%%%%%%%%%%%%%%
\multirow{3}{*}{20} & $h_0^5x_{127,15}$ & $d_{2}^{-1}$ & $h_0^4x_{128,14}$ \\\cline{2-4}
 & $d_0e_0[\Delta\Delta_1g]$ & $d_{4}$ & $d_0Pd_0M^2$ \\\cline{2-4}
 & $h_0^{19}h_7$ & $d_{3}$ & $h_0^5x_{126,18}$ \\\hline\hline
\multirow{5}{*}{19} & $h_0^4x_{127,15}$ & $d_{2}^{-1}$ & $h_0^3x_{128,14}$ \\\cline{2-4}
 & $h_1x_{126,18}$ & & Permanent \\\cline{2-4}
 & $e_0x_{110,15}$ & $d_{4}$ & $x_{126,23}$ \\\cline{2-4}
 & $h_1x_{126,18,2}$ & $d_{3}$ & $h_0x_{126,21}+h_0^4x_{126,18}$ \\\cline{2-4}
 & $h_0^{18}h_7$ & $d_{3}$ & $h_0^4x_{126,18}$ \\\hline\hline
\multirow{5}{*}{18} & $h_0^3x_{127,15}$ & $d_{2}^{-1}$ & $h_0^2x_{128,14}$ \\\cline{2-4}
 & $h_0^3h_6x_{64,14}$ & $d_{2}^{-1}$ & $h_0^2h_6x_{65,13}$ \\\cline{2-4}
 & $g^2\Delta h_1H_1$ & $d_{3}^{-1}$ & $gx_{108,11}$ \\\cline{2-4}
 & $h_1x_{126,17}$ & & Permanent \\\cline{2-4}
 & $h_0^{17}h_7$ & $d_{3}$ & $h_0^3x_{126,18}$ \\\hline\hline
\multirow{5}{*}{17} & $h_0^2h_2x_{124,14}$ & $d_{2}^{-1}$ & $x_{128,15}$ \\\cline{2-4}
 & $h_0^2x_{127,15}$ & $d_{2}^{-1}$ & $x_{128,15}+h_0x_{128,14}$ \\\cline{2-4}
 & $h_0^2h_6x_{64,14}$ & $d_{2}^{-1}$ & $h_0h_6x_{65,13}$ \\\cline{2-4}
 & $gx_{107,13}$ & $d_{4}^{-1}$ & $h_0h_3x_{121,11}$ \\\cline{2-4}
 & $h_0^{16}h_7$ & $d_{3}$ & $h_0^2x_{126,18}$ \\\hline\hline
\multirow{6}{*}{16} & $h_0x_{127,15}+h_0h_2x_{124,14}$ & $d_{2}^{-1}$ & $x_{128,14}$ \\\cline{2-4}
 & $h_0h_6x_{64,14}$ & $d_{2}^{-1}$ & $h_6x_{65,13}$ \\\cline{2-4}
 & $h_0h_2x_{124,14}$ & & Permanent \\\cline{2-4}
 & $x_{127,16}$ & & Permanent \\\cline{2-4}
 & $h_0^{15}h_7$ & $d_{3}$ & $h_0x_{126,18}$ \\\cline{2-4}
 & $gx_{107,12}$ & $d_{3}$ & $g^3x_{66,7}$ \\\hline\hline
\multirow{5}{*}{15} & $h_1x_{126,14}$ & $d_{2}^{-1}$ & $x_{128,13,2}$ \\\cline{2-4}
 & $h_2x_{124,14}$ & & Permanent \\\cline{2-4}
 & $x_{127,15}$ & $d_{5}$ & $d_0x_{112,16}$ \\\cline{2-4}
 & $h_0^{14}h_7$ & $d_{2}$ & $h_0^{15}h_6^2$ \\\cline{2-4}
 & $h_6x_{64,14}$ & $d_{2}$ & $h_1^2x_{124,15}$ \\\hline\hline
\multirow{3}{*}{14} & $h_0g\Delta h_6g$ & $d_{2}^{-1}$ & $x_{128,12,2}$ \\\cline{2-4}
 & $h_0h_3x_{120,12}$ & $d_{2}^{-1}$ & $h_3x_{121,11}$ \\\cline{2-4}
 & $h_0^{13}h_7$ & $d_{2}$ & $h_0^{14}h_6^2$ \\\hline\hline
\multirow{5}{*}{13} & $h_0^3x_{127,10}$ & $d_{2}^{-1}$ & $h_0x_{128,10}$ \\\cline{2-4}
 & $g\Delta h_6g$ & $d_{3}^{-1}$ & $x_{128,10,2}$ \\\cline{2-4}
 & $h_3x_{120,12}$ & $d_{4}^{-1}$ & $h_2x_{125,8,2}$ \\\cline{2-4}
 & $x_{127,13}$ & $d_{3}$ & $h_0^2D_2x_{68,8}$ \\\cline{2-4}
 & $h_0^{12}h_7$ & $d_{2}$ & $h_0^{13}h_6^2$ \\\hline\hline
\multirow{3}{*}{12} & $h_0^2x_{127,10}$ & $d_{2}^{-1}$ & $x_{128,10}$ \\\cline{2-4}
 & $h_1x_{126,11}$ & $d_{3}^{-1}$ & $h_1x_{127,8}$ \\\cline{2-4}
 & $h_0^{11}h_7$ & $d_{2}$ & $h_0^{12}h_6^2$ \\\hline\hline
\multirow{4}{*}{11} & $h_0h_3x_{120,9}$ & $d_{3}^{-1}$ & $h_3D_2h_6$ \\\cline{2-4}
 & $h_0h_2x_{124,9}$ & & Permanent \\\cline{2-4}
 & $h_0x_{127,10}$ & & Permanent \\\cline{2-4}
 & $h_0^{10}h_7$ & $d_{2}$ & $h_0^{11}h_6^2$ \\\hline\hline
\multirow{6}{*}{10} & $h_1^2x_{125,8}$ & & Permanent \\\cline{2-4}
 & $h_2x_{124,9}+h_0^2x_{127,8}$ & $d_{6}$ & $h_1^2x_{124,14}$ \\\cline{2-4}
 & $h_0^2x_{127,8}$ & $d_{4}$ & $x_{126,14}$ \\\cline{2-4}
 & $x_{127,10}$ & $d_{2}$ & $d_1x_{94,8}$ \\\cline{2-4}
 & $h_3x_{120,9}$ & $d_{2}$ & $h_0x_{126,11}$ \\\cline{2-4}
 & $h_0^9h_7$ & $d_{2}$ & $h_0^{10}h_6^2$ \\\hline
  \end{tabular}}
\caption{The classical Adams spectral sequence of $S^0$ for $10 \le s \le 20$ in stem 127}
\label{Table:S127.20}
\end{table}

\begin{table}
    \centering
\scalebox{0.7}{\begin{tabular}{|c|l|l|l|}\hline %%%%%%%%%%%%%%%%%%%%%%
\multirow{5}{*}{9} & $h_0^2x_{127,7}$ & $d_{2}^{-1}$ & $h_0x_{128,6}$ \\\cline{2-4}
 & $h_1x_{126,8}$ & & Permanent \\\cline{2-4}
 & $h_1x_{126,8,2}$ & $d_{5}$ & $h_1h_3x_{118,12}$ \\\cline{2-4}
 & $h_0x_{127,8}$ & $d_{2}$ & $h_0^2x_{126,9}$ \\\cline{2-4}
 & $h_0^8h_7$ & $d_{2}$ & $h_0^9h_6^2$ \\\hline\hline
\multirow{7}{*}{8} & $h_0x_{127,7,2}+h_0x_{127,7}+h_0^2x_{127,6}$ & $d_{2}^{-1}$ & $x_{128,6}$ \\\cline{2-4}
 & $h_0^2x_{127,6}$ & $d_{2}^{-1}$ & $h_0x_{128,5}$ \\\cline{2-4}
 & $h_2h_6A$ & & Permanent \\\cline{2-4}
 & $h_2x_{124,7}$ & $d_{9}$ & $?$ \\\cline{2-4}
 & $x_{127,8}$ & $d_{2}$ & $h_0x_{126,9}$ \\\cline{2-4}
 & $h_0x_{127,7}$ & $d_{2}$ & $h_0^2x_{126,8}$ \\\cline{2-4}
 & $h_0^7h_7$ & $d_{2}$ & $h_0^8h_6^2$ \\\hline\hline
\multirow{5}{*}{7} & $h_0x_{127,6}$ & $d_{2}^{-1}$ & $x_{128,5}$ \\\cline{2-4}
 & $h_1x_{126,6}$ & $d_{10}$ & $?$ \\\cline{2-4}
 & $x_{127,7,2}+x_{127,7}$ & $d_{3}$ & $?$ \\\cline{2-4}
 & $x_{127,7}$ & $d_{2}$ & $h_0x_{126,8}$ \\\cline{2-4}
 & $h_0^6h_7$ & $d_{2}$ & $h_0^7h_6^2$ \\\hline\hline
\multirow{2}{*}{6} & $x_{127,6}$ & $d_{4}$ & $?$ \\\cline{2-4}
 & $h_0^5h_7$ & $d_{2}$ & $h_0^6h_6^2$ \\\hline\hline
\multirow{1}{*}{5} & $h_0^4h_7$ & $d_{2}$ & $h_0^5h_6^2$ \\\hline\hline
\multirow{1}{*}{4} & $h_0^3h_7$ & $d_{2}$ & $h_0^4h_6^2$ \\\hline\hline
\multirow{2}{*}{3} & $h_1h_6^2$ & $d_{14}$ & $?$ \\\cline{2-4}
 & $h_0^2h_7$ & $d_{2}$ & $h_0^3h_6^2$ \\\hline\hline
\multirow{1}{*}{2} & $h_0h_7$ & $d_{2}$ & $h_0^2h_6^2$ \\\hline\hline
\multirow{1}{*}{1} & $h_7$ & $d_{2}$ & $h_0h_6^2$ \\\hline\hline
\multirow{1}{*}{0} & \multicolumn{3}{c|}{}\\\hline
 \end{tabular}}
\caption{The classical Adams spectral sequence of $S^0$ for $s \le 9$ in stem 127}
\label{Table:S127.9}
\end{table}

\begin{table}
    \centering
\scalebox{0.75}{\begin{tabular}{|c|l|l|l|}\hline 
$s$ & Elements & $d_r$ & value\\\hline\hline
\multirow{4}{*}{14} & $h_0^{12}h_6^2[0]$ & $d_{2}^{-1}$ & $h_0^{11}h_7[0]$ \\\cline{2-4}
 & $x_{126,14}[0]$ & $d_{3}^{-1}$ & $(((x_{123,11,2})+(x_{123,11})+h_0 h_6 [B_4])[4])$ \\\cline{2-4}
 & $Q_2D_2(h_3[4])+D_2x_{68,8}[0]$ & $d_{3}^{-1}$ & $(((x_{123,11,2})+h_5 (x_{92,10}))[4])$ \\\cline{2-4}
 & $D_2x_{68,8}[0]$ & $d_{2}$ & $h_0Q_2x_{68,8}[0]$ \\\hline\hline
\multirow{6}{*}{13} & $h_0^{11}h_6^2[0]$ & $d_{2}^{-1}$ & $h_0^{10}h_7[0]$ \\\cline{2-4}
 & $h_1x_{120,11}(h_1[4])$ & $d_{12}$ & $?$ \\\cline{2-4}
 & $h_1x_{125,12,2}[0]$ & $d_{5}$ & $d_0^2x_{97,10}[0]$ \\\cline{2-4}
 & $h_6x_{56,10}(h_0 h_2[4])$ & $d_{3}$ & $h_1x_{124,15}[0]$ \\\cline{2-4}
 & $(((x_{122,13})+h_1^2 (x_{120,11})+h_0^2 h_6 (Md_0))[4])$ & $d_{2}$ & $x_{125,15}[0]+h_0^3x_{125,12}[0]$ \\\cline{2-4}
 & $h_0h_3x_{119,11}[0]$ & $d_{2}$ & $h_0^3x_{125,12}[0]$ \\\hline\hline
\multirow{4}{*}{12} & $d_1x_{94,8}[0]$ & $d_{2}^{-1}$ & $x_{127,10}[0]$ \\\cline{2-4}
 & $h_0^{10}h_6^2[0]$ & $d_{2}^{-1}$ & $h_0^9h_7[0]$ \\\cline{2-4}
 & $((h_5 (x_{91,11})+h_0 (x_{122,11}))[4])$ & $d_{3}$ & $Q_2x_{68,8}[0]$ \\\cline{2-4}
 & $h_3x_{119,11}[0]$ & $d_{2}$ & $h_0^2x_{125,12}[0]$ \\\hline\hline
\multirow{4}{*}{11} & $h_0^9h_6^2[0]$ & $d_{2}^{-1}$ & $h_0^8h_7[0]$ \\\cline{2-4}
 & $x_{126,11}[0]$ & $d_{3}^{-1}$ & $x_{127,8}[0]$ \\\cline{2-4}
 & $h_1x_{125,10,2}[0]+h_1x_{125,10}[0]$ & & Permanent \\\cline{2-4}
 & $h_1x_{125,10}[0]$ & $d_{14}$ & $?$ \\\hline\hline
\multirow{3}{*}{10} & $h_0^2x_{126,8}[0]$ & $d_{2}^{-1}$ & $h_0x_{127,7}[0]$ \\\cline{2-4}
 & $h_0^8h_6^2[0]$ & $d_{2}^{-1}$ & $h_0^7h_7[0]$ \\\cline{2-4}
 & $x_{126,10}[0]$ & $d_{3}$ & $nx_{94,8}[0]$ \\\hline\hline
\multirow{5}{*}{9} & $h_0x_{126,8}[0]$ & $d_{2}^{-1}$ & $x_{127,7}[0]$ \\\cline{2-4}
 & $h_0^7h_6^2[0]$ & $d_{2}^{-1}$ & $h_0^6h_7[0]$ \\\cline{2-4}
 & $h_1x_{125,8}[0]$ & $d_{16}$ & $?$ \\\cline{2-4}
 & $h_0x_{126,8,3}[0]$ & $d_{4}$ & $h_0x_{125,12}[0]$ \\\cline{2-4}
 & $x_{126,9}[0]$ & $d_{3}$ & $h_0^4x_{125,8}[0]$ \\\hline %\hline
%\multicolumn{4}{|c|}{\text{stem}=126} \\\hline
%\multicolumn{4}{|c|}{\TitleCellColor $E_2^{*,*}(S^0/\nu)$}\\\hline
\end{tabular}}
\caption{The classical Adams spectral sequence of $S^0/\nu$ for $9 \le s \le 14$ in stem 126}
\label{Table:Cnu126}
\end{table}

\makebibliography
\end{document}

%% file: cw_data.tex
\begin{table}
\begin{tabular}{|l|c|}\hline
CW spectra & max(t)\\\hline
\texttt{S0} & 261\\\hline
\texttt{tmf} & 261\\\hline
\texttt{C2} & 200\\\hline
\texttt{Ceta} & 200\\\hline
\texttt{Cnu} & 200\\\hline
\texttt{Csigma} & 200\\\hline
\texttt{CW\_2\_eta} & 200\\\hline
\texttt{CW\_eta\_2} & 200\\\hline
\texttt{CW\_eta\_nu} & 200\\\hline
\texttt{CW\_nu\_eta} & 200\\\hline
\texttt{CW\_sigma\_nu} & 200\\\hline
\texttt{CW\_nu\_sigma} & 200\\\hline
\texttt{CW\_2\_eta\_nu} & 200\\\hline
\texttt{CW\_nu\_eta\_2} & 200\\\hline
\texttt{CW\_sigma\_nu\_eta} & 200\\\hline
\texttt{CW\_eta\_nu\_sigma} & 200\\\hline
\texttt{CW\_sigma\_nu\_eta\_2} & 200\\\hline
\texttt{CW\_2\_eta\_nu\_sigma} & 200\\\hline
\texttt{Csigmasq} & 200\\\hline
\texttt{C2h4} & 200\\\hline
\texttt{DC2h4} & 198\\\hline
\texttt{Ctheta4} & 200\\\hline
\texttt{C2h5} & 200\\\hline
\texttt{DC2h5} & 200\\\hline
\texttt{Ctheta5} & 200\\\hline
\end{tabular}
\begin{tabular}{|l|c|}\hline
CW spectra & max(t)\\\hline
\texttt{C2h6} & 200\\\hline
\texttt{DC2h6} & 200\\\hline
\texttt{C2\_C2} & 180\\\hline
\texttt{Ceta\_Ceta} & 160\\\hline
\texttt{Cnu\_Cnu} & 180\\\hline
\texttt{Csigma\_Csigma} & 180\\\hline
\texttt{CW\_2sigma\_sigma} & 180\\\hline
\texttt{CW\_sigma\_2sigma} & 180\\\hline
\texttt{C2sigma} & 200\\\hline
\texttt{CW\_2\_V\_eta} & 200\\\hline
\texttt{CW\_2\_A\_eta} & 180\\\hline
\texttt{Joker} & 200\\\hline
\texttt{CW\_eta\_2\_eta\_Eq\_2\_nu} & 200\\\hline
\texttt{CW\_eta\_2\_eta\_Eq\_nu\_2} & 200\\\hline
\texttt{C2\_Ceta} & 200\\\hline
\texttt{RP3\_6} & 199\\\hline
\texttt{Fphi} & 200\\\hline
\texttt{RP1\_4} & 200\\\hline
\texttt{RP1\_6} & 183\\\hline
\texttt{RP1\_8} & 180\\\hline
\texttt{RP1\_10} & 180\\\hline
\texttt{RP1\_12} & 178\\\hline
\texttt{RP1\_256} & 200\\\hline
\texttt{RP3\_256} & 160\\\hline
 & \\\hline
\end{tabular}
\caption{CW spectra}\label{tb:CW1}
\end{table}

%% file: cw_map_data.tex
\begin{table}\begin{tabular}{|l|c|}\hline
map & max(t)\\\hline
\texttt{S0\_\_tmf} & 261\\\hline
\texttt{C2\_\_S0} & 200\\\hline
\texttt{Ceta\_\_S0} & 200\\\hline
\texttt{Cnu\_\_S0} & 200\\\hline
\texttt{Csigma\_\_S0} & 200\\\hline
\texttt{C2\_\_CW\_2\_eta} & 180\\\hline
\texttt{Ceta\_\_CW\_eta\_2} & 180\\\hline
\texttt{Ceta\_\_CW\_eta\_nu} & 180\\\hline
\texttt{Cnu\_\_CW\_nu\_eta} & 180\\\hline
\texttt{Cnu\_\_CW\_nu\_sigma} & 200\\\hline
\texttt{CW\_nu\_sigma\_\_Csigma} & 199\\\hline
\texttt{CW\_sigma\_nu\_\_Cnu} & 199\\\hline
\texttt{CW\_2\_eta\_\_Ceta} & 200\\\hline
\texttt{CW\_eta\_2\_\_C2} & 180\\\hline
\texttt{CW\_eta\_nu\_\_Cnu} & 200\\\hline
\texttt{CW\_nu\_eta\_\_Ceta} & 200\\\hline
\texttt{C2\_\_Q\_CW\_nu\_eta\_2} & 194\\\hline
\texttt{Cnu\_\_Q\_CW\_2\_eta\_nu} & 197\\\hline
\texttt{CW\_2\_eta\_nu\_\_CW\_eta\_nu} & 190\\\hline
\texttt{CW\_2\_eta\_\_CW\_2\_eta\_nu} & 200\\\hline
\texttt{CW\_eta\_2\_\_Q\_CW\_nu\_eta\_2} & 200\\\hline
\texttt{CW\_eta\_nu\_\_Q\_CW\_2\_eta\_nu} & 200\\\hline
\texttt{CW\_nu\_eta\_2\_\_CW\_eta\_2} & 200\\\hline
\texttt{CW\_nu\_eta\_\_CW\_nu\_eta\_2} & 187\\\hline
\texttt{C2\_\_CW\_sigma\_nu\_by\_eta} & 179\\\hline
\texttt{Ceta\_\_Csigma\_by\_nu} & 188\\\hline
\texttt{Csigma\_\_Ceta\_by\_nu} & 194\\\hline
\texttt{Csigma\_\_CW\_2\_eta\_by\_nu} & 193\\\hline
\texttt{CW\_2\_eta\_nu\_sigma\_\_CW\_eta\_nu\_sigma} & 200\\\hline
\texttt{CW\_2\_eta\_nu\_\_CW\_2\_eta\_nu\_sigma} & 200\\\hline
\texttt{CW\_eta\_2\_\_Csigma\_by\_nu} & 188\\\hline
\texttt{CW\_eta\_nu\_sigma\_\_CW\_nu\_sigma} & 200\\\hline
\texttt{CW\_eta\_nu\_sigma\_\_S0\_by\_2} & 200\\\hline
\texttt{CW\_eta\_nu\_\_CW\_eta\_nu\_sigma} & 200\\\hline
\texttt{CW\_nu\_eta\_2\_\_S0\_by\_sigma} & 200\\\hline
\texttt{CW\_nu\_eta\_\_S0\_by\_sigma} & 197\\\hline
\texttt{CW\_nu\_sigma\_\_C2\_by\_eta} & 197\\\hline
\texttt{CW\_nu\_sigma\_\_S0\_by\_eta} & 199\\\hline
\texttt{CW\_sigma\_nu\_eta\_2\_\_CW\_nu\_eta\_2} & 200\\\hline
\texttt{CW\_sigma\_nu\_eta\_\_CW\_nu\_eta} & 200\\\hline
\texttt{CW\_sigma\_nu\_eta\_\_CW\_sigma\_nu\_eta\_2} & 200\\\hline
\texttt{Ceta\_\_Q\_CW\_2\_eta} & 200\\\hline
\texttt{Ceta\_\_Q\_CW\_nu\_eta} & 200\\\hline
\texttt{Cnu\_\_Q\_CW\_eta\_nu} & 200\\\hline
\texttt{Cnu\_\_Q\_CW\_sigma\_nu} & 200\\\hline
\end{tabular}\caption{Maps between CW spectra}\label{tb:map1}
\end{table}

\begin{table}\begin{tabular}{|l|c|}\hline
map & max(t)\\\hline
\texttt{Csigma\_\_CW\_sigma\_nu} & 190\\\hline
\texttt{Csigma\_\_Q\_CW\_nu\_sigma} & 200\\\hline
\texttt{CW\_sigma\_nu\_\_CW\_sigma\_nu\_eta} & 200\\\hline
\texttt{Csigma\_\_Csigmasq} & 193\\\hline
\texttt{Csigmasq\_\_S0} & 200\\\hline
\texttt{C2h4\_\_Csigmasq} & 200\\\hline
\texttt{Csigmasq\_\_DC2h4} & 196\\\hline
\texttt{Csigmasq\_\_Csigma} & 200\\\hline
\texttt{Csigmasq\_\_C2\_by\_sigmasq} & 184\\\hline
\texttt{Csigmasq\_\_Q\_C2h4} & 200\\\hline
\texttt{Csigmasq\_\_S0\_by\_sigmasq} & 200\\\hline
\texttt{C2\_\_C2h4} & 200\\\hline
\texttt{DC2h4\_\_C2} & 193\\\hline
\texttt{C2h4\_\_S0} & 200\\\hline
\texttt{C2\_\_Q\_DC2h4} & 200\\\hline
\texttt{C2h4\_\_S0\_by\_sigmasq} & 199\\\hline
\texttt{C2\_\_C2h5} & 200\\\hline
\texttt{C2h5\_\_Ctheta4} & 200\\\hline
\texttt{Ctheta4\_\_DC2h5} & 200\\\hline
\texttt{DC2h5\_\_C2} & 200\\\hline
\texttt{C2h5\_\_S0} & 200\\\hline
\texttt{C2\_\_Q\_DC2h5} & 200\\\hline
\texttt{C2h5\_\_S0\_by\_theta4} & 200\\\hline
\texttt{Ctheta4\_\_Q\_C2h5} & 200\\\hline
\texttt{Ctheta4\_\_S0} & 200\\\hline
\texttt{Ctheta4\_\_C2\_by\_theta4} & 168\\\hline
\texttt{Ctheta4\_\_S0\_by\_theta4} & 200\\\hline
\texttt{C2\_\_C2h6} & 200\\\hline
\texttt{C2h6\_\_Ctheta5} & 200\\\hline
\texttt{Ctheta5\_\_DC2h6} & 200\\\hline
\texttt{DC2h6\_\_C2} & 200\\\hline
\texttt{C2h6\_\_S0} & 200\\\hline
\texttt{C2\_\_Q\_DC2h6} & 198\\\hline
\texttt{C2h6\_\_S0\_by\_theta5} & 198\\\hline
\texttt{Ctheta5\_\_Q\_C2h6} & 200\\\hline
\texttt{Ctheta5\_\_S0} & 200\\\hline
\texttt{C2\_\_Q\_CW\_2\_theta5\_2\_Eq\_eta\_theta5} & 136\\\hline
\texttt{Ctheta5\_\_Q\_CW\_2\_theta5\_2\_Eq\_eta\_theta5} & 198\\\hline
\texttt{DC2h6\_\_Q\_CW\_2\_theta5\_2\_Eq\_eta\_theta5} & 191\\\hline
\texttt{C2\_C2\_\_C2} & 180\\\hline
\texttt{C2\_\_C2\_C2} & 180\\\hline
\texttt{CW\_2\_eta\_\_Q\_RP1\_6} & 180\\\hline
\texttt{Ceta\_Ceta\_\_Ceta} & 157\\\hline
\texttt{Ceta\_\_Ceta\_Ceta} & 160\\\hline
\texttt{Cnu\_Cnu\_\_Cnu} & 180\\\hline
\end{tabular}\caption{Maps between CW spectra}
\end{table}

\begin{table}\begin{tabular}{|l|c|}\hline
map & max(t)\\\hline
\texttt{Cnu\_\_Cnu\_Cnu} & 180\\\hline
\texttt{Csigma\_Csigma\_\_CW\_sigma\_2sigma} & 180\\\hline
\texttt{Csigma\_Csigma\_\_Csigma} & 180\\\hline
\texttt{Csigma\_\_CW\_sigma\_2sigma} & 180\\\hline
\texttt{Csigma\_\_Csigma\_Csigma} & 180\\\hline
\texttt{C2sigma\_\_S0} & 200\\\hline
\texttt{C2sigma\_\_CW\_2sigma\_sigma} & 160\\\hline
\texttt{CW\_sigma\_2sigma\_\_C2sigma} & 180\\\hline
\texttt{C2\_C2\_\_CW\_2\_V\_eta} & 160\\\hline
\texttt{C2\_\_Q\_CW\_2\_A\_eta} & 139\\\hline
\texttt{CW\_2\_A\_eta\_\_C2\_C2} & 160\\\hline
\texttt{Ceta\_\_Q\_CW\_2\_A\_eta} & 160\\\hline
\texttt{Ceta\_\_Q\_CW\_eta\_A\_nu} & 160\\\hline
\texttt{Cnu\_\_Q\_CW\_eta\_A\_nu} & 160\\\hline
\texttt{Cnu\_\_Q\_CW\_nu\_A\_sigma} & 160\\\hline
\texttt{Csigma\_\_Q\_CW\_nu\_A\_sigma} & 160\\\hline
\texttt{CW\_2\_V\_eta\_\_S1} & 160\\\hline
\texttt{CW\_2\_V\_eta\_\_S2} & 160\\\hline
\texttt{Ceta\_\_Joker} & 193\\\hline
\texttt{Ceta\_\_Q\_Joker} & 198\\\hline
\texttt{CW\_2\_eta\_\_Joker} & 200\\\hline
\texttt{CW\_eta\_2\_\_Q\_Joker} & 199\\\hline
\texttt{Joker\_\_Ceta} & 200\\\hline
\texttt{Joker\_\_CW\_eta\_2} & 200\\\hline
\texttt{C2\_\_Q\_Joker} & 197\\\hline
\texttt{CW\_2\_V\_eta\_\_Joker} & 192\\\hline
\texttt{Joker\_\_CW\_2\_A\_eta} & 182\\\hline
\texttt{CW\_2\_A\_eta\_\_Q\_Joker} & 180\\\hline
\texttt{CW\_sigma\_nu\_eta\_2\_\_tmf} & 200\\\hline
\texttt{C2\_Ceta\_\_C2} & 200\\\hline
\texttt{C2\_Ceta\_\_Ceta} & 200\\\hline
\texttt{C2\_\_C2\_Ceta} & 180\\\hline
\texttt{Ceta\_\_C2\_Ceta} & 180\\\hline
\texttt{Cnu\_\_Q\_CW\_eta\_nu\_eta\_Eq\_nu\_nu} & 196\\\hline
\texttt{Csigma\_\_Q\_CW\_nu\_sigma\_nu\_Eq\_sigma\_sigma} & 192\\\hline
\texttt{CW\_nu\_eta\_\_Q\_CW\_eta\_nu\_eta\_Eq\_nu\_nu} & 195\\\hline
\texttt{CW\_sigma\_nu\_\_Q\_CW\_nu\_sigma\_nu\_Eq\_sigma\_sigma} & 194\\\hline
\texttt{C2\_\_CW\_eta\_2\_eta\_Eq\_2\_nu} & 193\\\hline
\texttt{C2\_\_Q\_CW\_eta\_2\_eta\_Eq\_nu\_2} & 196\\\hline
\texttt{Cnu\_\_Q\_CW\_eta\_2\_eta\_Eq\_2\_nu} & 199\\\hline
\texttt{CW\_2\_eta\_\_Q\_CW\_eta\_2\_eta\_Eq\_2\_nu} & 198\\\hline
\texttt{CW\_2\_eta\_\_Q\_CW\_eta\_2\_eta\_Eq\_nu\_2} & 198\\\hline
\texttt{CW\_eta\_2\_eta\_Eq\_2\_nu\_\_Cnu} & 193\\\hline
\texttt{CW\_eta\_2\_eta\_Eq\_2\_nu\_\_CW\_2\_eta} & 193\\\hline
\texttt{CW\_eta\_2\_eta\_Eq\_nu\_2\_\_C2} & 200\\\hline
\end{tabular}\caption{Maps between CW spectra}
\end{table}

\begin{table}\begin{tabular}{|l|c|}\hline
map & max(t)\\\hline
\texttt{CW\_eta\_2\_eta\_Eq\_nu\_2\_\_CW\_2\_eta} & 200\\\hline
\texttt{CW\_eta\_2\_\_CW\_eta\_2\_eta\_Eq\_2\_nu} & 193\\\hline
\texttt{CW\_eta\_2\_\_CW\_eta\_2\_eta\_Eq\_nu\_2} & 200\\\hline
\texttt{Cnu\_\_CW\_eta\_2\_eta\_Eq\_nu\_2} & 200\\\hline
\texttt{C2\_\_CW\_2\_V\_eta} & 177\\\hline
\texttt{CW\_2\_V\_eta\_\_C2\_Ceta} & 200\\\hline
\texttt{CW\_2\_V\_eta\_\_CW\_eta\_2\_eta\_Eq\_2\_nu} & 200\\\hline
\texttt{Ceta\_\_CW\_2\_V\_eta} & 200\\\hline
\texttt{C2\_Ceta\_\_CW\_2\_A\_eta} & 181\\\hline
\texttt{CW\_2\_A\_eta\_\_C2} & 180\\\hline
\texttt{CW\_2\_A\_eta\_\_Ceta} & 180\\\hline
\texttt{CW\_2\_A\_eta\_\_Q\_C2\_Ceta} & 180\\\hline
\texttt{CW\_2\_A\_eta\_\_Q\_CW\_eta\_2\_eta\_Eq\_nu\_2} & 180\\\hline
\texttt{CW\_eta\_2\_eta\_Eq\_nu\_2\_\_CW\_2\_A\_eta} & 183\\\hline
\texttt{Ceta\_\_Q\_CW\_eta\_2\_eta\_Eq\_2\_nu} & 197\\\hline
\texttt{CW\_eta\_2\_\_Fphi} & 199\\\hline
\texttt{RP5\_8\_\_RP3\_4} & 173\\\hline
\texttt{CW\_eta\_2\_\_RP1\_6} & 180\\\hline
\texttt{RP1\_6\_\_CW\_2\_eta} & 180\\\hline
\texttt{Fphi\_\_RP1\_256} & 196\\\hline
\texttt{RP1\_256\_\_S0} & 199\\\hline
\texttt{RP1\_4\_\_CP1\_2} & 200\\\hline
\texttt{RP1\_256\_\_RPm7\_0} & 172\\\hline
\texttt{RP1\_10\_\_RP7\_10} & 180\\\hline
\texttt{RP9\_16\_\_RP1\_8} & 156\\\hline
\texttt{RP3\_4\_\_RP1\_2} & 198\\\hline
\texttt{RP5\_6\_\_RP3\_4} & 198\\\hline
\texttt{RP1\_4\_\_RP2\_4} & 160\\\hline
\texttt{RP3\_5\_\_RP3\_6} & 160\\\hline
\texttt{RP3\_4\_\_RP3\_6} & 188\\\hline
\texttt{RP1\_2\_\_RP1\_4} & 199\\\hline
\texttt{RP1\_4\_\_RP1\_6} & 183\\\hline
\texttt{RP1\_6\_\_RP1\_8} & 175\\\hline
\texttt{RP1\_8\_\_RP1\_10} & 174\\\hline
\texttt{RP1\_10\_\_RP1\_12} & 178\\\hline
\texttt{RP1\_12\_\_RP1\_256} & 178\\\hline
\texttt{RP1\_6\_\_RP3\_6} & 179\\\hline
\texttt{RP3\_6\_\_RP5\_6} & 191\\\hline
\texttt{RP1\_4\_\_RP3\_4} & 200\\\hline
\texttt{RP1\_6\_\_RP5\_6} & 162\\\hline
\texttt{RP1\_8\_\_RP5\_8} & 176\\\hline
\texttt{RP1\_12\_\_RP9\_12} & 178\\\hline
\texttt{RP1\_256\_\_RP3\_256} & 160\\\hline
\texttt{RP3\_6\_\_RP3\_256} & 140\\\hline
\texttt{RP3\_256\_\_RP1\_2} & 160\\\hline
\end{tabular}\caption{Maps between CW spectra}\label{tb:map2}
\end{table}